\documentclass[12pt,reqno]{amsart} 
\usepackage{times}

\newcommand{\Char}{\operatorname{char}}
\newcommand{\codim}{\operatorname{codim}}
\newcommand{\Dbar}{\overline{D}}
\newcommand{\Dibar}{\overline{D_1}}
\newcommand{\Diibar}{\overline{D_2}}
\newcommand{\Gal}{\operatorname{Gal}}
\newcommand{\GL}{\operatorname{GL}}
\newcommand{\PGL}{\operatorname{PGL}}
\newcommand{\Kbar}{\overline{K}}
\newcommand{\la}{\lambda}
\newcommand{\lcm}{{\operatorname{lcm}}}
\newcommand{\LL}{{\mathcal{L}}}
\newcommand{\Mat}{\operatorname{Mat}}
\newcommand{\one}{{\mathbf{1}}}
\newcommand{\pf}{\operatorname{pf}}
\newcommand{\PP}{{\mathbb P}}
\newcommand{\R}{{\mathbb R}}
\newcommand{\rank}{\operatorname{rank}}
\newcommand{\ra}{{\longrightarrow}}
\newcommand{\Sec}{\operatorname{Sec}}
\newcommand{\sign}{\operatorname{sign}}
\newcommand{\SL}{\operatorname{SL}}
\newcommand{\vP}{{\bf P}}
\newcommand{\vv}{{\bf v}}
\newcommand{\vQ}{{\bf Q}}
\newcommand{\Z}{{\mathbb Z}}

\newtheorem{Theorem}{Theorem}[section]
\newtheorem{Corollary}[Theorem]{Corollary}
\newtheorem{Lemma}[Theorem]{Lemma}
\newtheorem{Proposition}[Theorem]{Proposition}

\newenvironment{Proof}{\par\noindent{\bf Proof:}}%
                      {\hspace*{\fill}\nobreak$\Box$\par\medskip}
\newenvironment{ProofOf}[1]{\par\noindent{\bf Proof of #1:}}%
                      {\hspace*{\fill}\nobreak$\Box$\par\medskip}

\theoremstyle{definition}
\newtheorem{Remark}[Theorem]{Remark}

\title[The Jacobian of a genus one curve]%
{\large A formula for the Jacobian of \\ a genus one curve of arbitrary degree}

\author{Tom~Fisher}
\address{University of Cambridge,
           DPMMS, Centre for Mathematical Sciences,
           Wilberforce Road, Cambridge CB3 0WB, UK}
\email{T.A.Fisher@dpmms.cam.ac.uk}
\date{14th October 2015}  

\begin{document}

\begin{abstract}
  We extend the formulae of classical invariant theory for the
  Jacobian of a genus one curve of degree $n \le 4$ to curves of
  arbitrary degree. To do this, we associate to each genus one normal
  curve of degree $n$, an $n \times n$ alternating matrix of quadratic
  forms in $n$ variables, that represents the invariant
  differential. We then exhibit the invariants we need as homogeneous
  polynomials of degrees $4$ and $6$ in the coefficients of the
  entries of this matrix.
\end{abstract}

\maketitle

\renewcommand{\theenumi}{\roman{enumi}}

\section*{Introduction}

Let $C$ be a smooth curve of genus one defined over a field $K$.  Its
Jacobian is an elliptic curve $E$ defined over the same field
$K$. However it is only if $C$ has a $K$-rational point that $C$ and
$E$ are isomorphic over $K$. Starting with equations for $C$ we would
like to compute a Weierstrass equation for $E$.

Let $D$ be a $K$-rational divisor on $C$ of degree $n \ge 1$.  It is
natural to split into cases according to the value of $n$.  If $n=1$
then $C$ has a $K$-rational point, and our task is that of writing an
elliptic curve in Weierstrass form. If $n \ge 2$ then the complete
linear system $|D|$ defines a morphism $C \to \PP^{n-1}$.  Explicitly,
the map is given by $(f_1: \ldots :f_n)$ where $f_1, \ldots, f_n$ are
a basis for the Riemann-Roch space $\LL(D)$.  If $n=2$ then $C$ is a
double cover of $\PP^1$ and is given by an equation of the form $y^2 =
F(x_1,x_2)$ where $F$ is a binary quartic. In this case
Weil~\cite{WHermite}, \cite{WEuler} showed that the classical
invariants of the binary quartic $F$ give a formula for the Jacobian.

If $n \ge 3$ then the morphism $C \to \PP^{n-1}$ is an embedding.  The
image is a {\em genus one normal curve} of degree $n$.  The word {\em
  normal} refers to the fact $C$ is projectively normal (see for
example \cite[Proposition IV.1.2]{Hulek}), i.e. if $H$ is the divisor
of a hyperplane section then the natural map
\begin{equation}
\label{projnormal}
S^d \LL(H)  \to \LL(d H)
\end{equation} 
is surjective for all $d \ge 1$.  If $n = 3$ then $C \subset \PP^2$ is
a smooth plane cubic, say with equation $F(x_1,x_2,x_3) = 0$. The
invariants of a ternary cubic $F$ were computed by
Aronhold~\cite{Aronhold}, and again Weil (in the notes
to~\cite{WHermite} in his collected papers) showed that these give a
formula for the Jacobian.  If $n = 4$ then $C \subset \PP^3$ is the
complete intersection of two quadrics.  If we represent these quadrics
by $4 \times 4$ symmetric matrices $A$ and $B$, then $F(x_1,x_2) =
\det(A x_1 + B x_2)$ is a binary quartic. The invariants of this
binary quartic again give a formula for the Jacobian. For further
details of these formulae in the cases $n=2,3,4$ see \cite{Mc+},
\cite{ARVT} or \cite{g1inv}.

If $n=5$ then $C \subset \PP^4$ is no longer a complete intersection,
and indeed the homogeneous ideal is generated by $5$ quadrics.  The
Buchsbaum-Eisenbud structure theorem \cite{BE1}, \cite{BE2} shows that
these quadrics may be written as the $4 \times 4$ Pfaffians of a $5
\times 5$ alternating matrix of linear forms. The space of all such
matrices is a $50$-dimensional affine space, with a natural action of
$\GL_5 \times \GL_5$.  In \cite{g1inv} we computed generators for the
ring of invariants and showed that they again give a formula for the
Jacobian.  In fact the invariants are too large to write down as
explicit polynomials, so instead we gave a practical algorithm for
evaluating them (based in part on the case $n=5$ of
Proposition~\ref{lem:deriv1}).  More recently, B. Gross
\cite{BGross} gave a uniform description of the invariants in the
cases $n=2,3,4,5$, using results of Vinberg, although this does not
appear to give any way of evaluating the invariants in the case $n=5$.

In this paper we extend these formulae for the Jacobian to genus one
normal curves of arbitrary degree.

Let $C \subset \PP^{n-1}$ be a genus one normal curve of degree $n \ge
3$.  Since $C$ has genus one, the space of regular differentials on
$C$ has dimension $1$, say spanned by $\omega$. We call $\omega$ an
{\em invariant differential}, since geometrically it is invariant
under all translation maps.  There is a linear map
\begin{equation}
\label{firstomega}
 \wedge^2 \LL(H) \to \LL(2H) \,; \quad 
   f \wedge g \mapsto \frac{f dg - g df}{\omega}.  
\end{equation}
Since~\eqref{projnormal} is surjective for $d=2$, we may represent
this map by an $n \times n$ alternating matrix of quadratic forms in
$x_1, \ldots, x_n$.  This matrix $\Omega$ represents $\omega$ in the
sense that
\[ \omega = \frac{ x_j^2 d(x_i/x_j) }{\Omega_{ij}(x_1, \ldots,x_n)}
\qquad \text{ for all } i \not= j. \] However if $n \ge 4$ then there
are quadrics vanishing on $C \subset \PP^{n-1}$ and so this
description does not determine $\Omega$ uniquely. Nonetheless we show,
by proving \cite[Conjecture~7.4]{invenqI}, that there is a canonical
choice of $\Omega$.  We then define polynomials $c_4$ and $c_6$ of
degrees $4$ and $6$ in the coefficients of the entries of $\Omega$,
and show that the Jacobian has Weierstrass equation
\[ y^2 = x^3 - 27 c_4(\Omega) x - 54 c_6(\Omega). \]

These main results are stated in Section~\ref{sec:stat}. In the next
two sections we show that $c_4$ and $c_6$ are invariants for the
appropriate action of $\GL_n$, and that they reduce to the previously
known formulae for $n \le 5$.  At this point the proof of our results
for any given value of $n$ is a finite calculation. However finding a
proof that works for all $n$ is more challenging.

In Section~\ref{sec:outline} we show that if we can find a matrix
$\Omega$ satisfying some apparently weaker hypotheses, then it will
satisfy the properties claimed in Theorem~\ref{thm1}.  For the actual
construction of $\Omega$ in Section~\ref{sec:pf1} we reduce to the
case $C$ is an elliptic curve $E$ embedded in $\PP^{n-1}$ via the
complete linear system $|n.0_E|$. At first we specify $\Omega$ as a
linear map $\wedge^2 \LL(n.0_E) \to S^2 \LL(n.0_E)$, and use this in
Section~\ref{sec:checkhyp} to complete the proof of
Theorem~\ref{thm1}. Then in Section~\ref{sec:pf2} we 
make a specific choice of basis for
$\LL(n.0_E)$, so that $\Omega$ 
becomes an alternating matrix of quadratic forms.  We
compute this matrix explicitly and, in Section~\ref{sec:scale},
prove the formula for the Jacobian by computing $c_4(\Omega)$ and
$c_6(\Omega)$. Much of the work here is in checking that the
invariants $c_4$ and $c_6$ are scaled correctly for all~$n$.

The description of $\Omega$ in Theorem~\ref{thm1} involves higher
secant varieties.  We quote any general results we need about these as
required.  Proofs, or references to the literature, are then given in
Section~\ref{sec:hsec}.

In future work we plan to study the space of all matrices $\Omega$.
This appears to be defined by $d_1+d_2$ quadrics in $\PP^{N-1}$ where
$N = (n^2-1) (n^2-4)/4$ and
\begin{align*}
d_1 &= (n^2-1)(n^2-4)(n^2-9)/36, \\
d_2 &= (n^2-1)^2 (n^2-9)/9.
\end{align*}
The numbers $N$, $d_1$ and $d_2$ are dimensions of irreducible 
representations for $\GL_n$.
Moreover, as suggested by Manjul Bhargava, we expect that $d_2$ of the
quadrics can be explained by an associative law, similar to that used
in \cite[Section 4]{quintic}.

We work throughout over a field $K$ of characteristic $0$, although it
would in fact be sufficient that the characteristic is not too small
compared to $n$. Except at the end of Section~\ref{sec:stat}, where we
give the application to computing Jacobians, we will assume that $K$
is algebraically closed.  For a projective variety $X$ we write $I(X)$
for its homogeneous ideal, and $T_P X$ for the tangent space at $P \in
X$.

\section{Statement of results}
\label{sec:stat}

Let $C \subset \PP^{n-1}$ be a genus one normal curve of degree $n \ge
3$.  For any integer $r \ge 1$ the $r$th {\em higher secant variety}
$\Sec^r C$ is the Zariski closure of the locus of all $(r-1)$-planes
through $r$ points on $C$.  For example, if $r=1$ then $\Sec^1 C = C$.
The codimension of $\Sec^r C$ in $\PP^{n-1}$ is $\max(n-2r,0)$. So
according as $n$ is odd or even there is a higher secant variety of
codimension $1$ or $2$.  If $n=2r+1$ then $\Sec^r C$ is a hypersurface
of degree $n$, whereas if $n=2r+2$ then $\Sec^r C$ is the complete
intersection of two forms of degree $r+1$.  In Section~\ref{sec:hsec}
we give references for these facts about higher secant varieties, and
also explain how to compute equations for $\Sec^r C$ from equations
for $C$.

We give the polynomial ring $R = K[x_1,\ldots,x_n]$ its usual grading
by degree, say $R = \oplus_d R_d$, and write $R(d)$ for the graded
$R$-module with $e$th graded piece $R_{d+e}$. Maps between graded free
$R$-modules are required to have relative degree $0$, and are labelled
by the matrices of forms that represent them. Our first main result is
\begin{Theorem}
\label{thm1}
Let $C \subset \PP^{n-1}$ be a genus one normal curve of degree $n \ge
3$.
\begin{enumerate}
\item If $n$ is odd, say $n=2r+1$, and $\Sec^r C = \{ F = 0 \}$ then
  there is a minimal free resolution
  \[ 0 \ra R(-2n) \stackrel{\nabla^T}{\ra} R(-n-1)^n
  \stackrel{\Omega}{\ra} R(-n+1)^n \stackrel{\nabla}{\ra} R \] where
  $\Omega$ is an $n \times n$ alternating matrix of quadratic forms
  and
  \[ \nabla = \nabla(F) = \begin{pmatrix} \frac{\partial F}{\partial
      x_1} & \cdots & \frac{\partial F}{\partial
      x_n} \end{pmatrix}. \]
\item If $n$ is even, say $n=2r+2$, and $\Sec^r C = \{ F_1 = F_2 =
  0\}$ then there is a minimal free resolution
  \[ 0 \ra R(-n)^2 \stackrel{\nabla^T}{\ra} R(\tfrac{-n-2}{2})^n
  \stackrel{\Omega}{\ra} R(\tfrac{-n+2}{2})^n \stackrel{\nabla}{\ra}
  R^2 \] where $\Omega$ is an $n \times n$ alternating matrix of
  quadratic forms and
  \[ \nabla = \nabla(F_1,F_2) = \begin{pmatrix} \vspace{0.5ex}
    \frac{\partial F_1}{\partial x_1} & \cdots
    & \frac{\partial F_1}{\partial x_n} \\
    \frac{\partial F_2}{\partial x_1} & \cdots & \frac{\partial
      F_2}{\partial x_n}
\end{pmatrix}. \]
\end{enumerate}
\end{Theorem}

We remarked in \cite[Section 7]{invenqI} that Theorem~\ref{thm1}(i)
follows from the Buchsbaum-Eisenbud structure theorem for Gorenstein
ideals of codimension $3$. In this paper we give a different proof,
not only so that it runs in parallel with
our proof of
Theorem~\ref{thm1}(ii), but also because this is needed for the proof
of Theorem~\ref{thm2}.

If the matrix $\Omega$ exists then, by the uniqueness of minimal free
resolutions (see for example \cite[Section 20.1]{Eisenbud} or
\cite[Section 7]{Peeva}), it is uniquely determined up to
scalars. Moreover starting from equations for $\Sec^r C$ we can solve
for $\Omega$ by linear algebra. The details are very similar to those
in \cite[Section 4]{5desc}.

Let $\Omega=(\Omega_{ij})$ be as specified in Theorem~\ref{thm1}. We
put
\begin{equation}
\label{MN}
M_{ij} = \sum_{r,s=1}^n 
\frac{\partial \Omega_{ir}}{\partial x_s} 
\frac{\partial \Omega_{js}}{\partial x_r} 
\quad \text{ and } \quad
N_{ijk} = \sum_{r=1}^n \frac{\partial M_{ij}}{\partial x_r} \Omega_{rk}. 
\end{equation}
We then define 
\begin{equation}
\label{c4}
c_4(\Omega) =  \frac{3 (n-2)^2}{2^4 n \binom{n+3}{5}}
\sum_{i,j,r,s=1}^n \frac{ \partial^2 M_{ij}}{\partial x_r \partial x_s} 
\frac{ \partial^2 M_{rs}}{\partial x_i \partial x_j} 
\end{equation}
and 
\begin{equation}
\label{c6}
c_6(\Omega) =  \frac{-(n-2)^3}{2^6 n \binom{n+5}{7}}
\sum_{i,j,k,r,s,t=1}^n \frac{ \partial^3 N_{ijk}}{\partial x_r \partial x_s
  \partial x_t } 
\frac{\partial^3 N_{rst}}{\partial x_i \partial x_j \partial x_k}. 
\end{equation}
Let $C_1$ and $C_2$ be genus one curves with invariant differentials
$\omega_1$ and $\omega_2$.  An isomorphism $\gamma : (C_1,\omega_1)
\to (C_2,\omega_2)$ is an isomorphism of curves $\gamma : C_1 \to C_2$
with $\gamma^* \omega_2 = \omega_1$.

\begin{Theorem}
\label{thm2}
Let $C \subset \PP^{n-1}$ be a genus one normal curve of degree $n \ge
3$, and let $\Omega$ be an alternating matrix of quadratic forms as
specified in Theorem~\ref{thm1}. Then
\begin{enumerate}
\item There is an invariant differential $\omega$ on $C$ such that
\begin{equation*}
\omega = \frac{ x_j^2 d(x_i/x_j) }{\Omega_{ij}(x_1, \ldots,x_n)} \qquad
\text{ for all } i \not= j. 
\end{equation*}
\item The pair $(C,\omega)$ is isomorphic (over $K = \Kbar$) to
\[ (y^2 = x^3 - 27 c_4(\Omega) x - 54 c_6(\Omega), 3 dx/y ). \]
\end{enumerate}
\end{Theorem}

The following corollary gives the application of Theorem~\ref{thm2} to
computing Jacobians. For this result only we drop our assumption that
$K$ is algebraically closed.

\begin{Corollary}
\label{compjac}
Let $C \subset \PP^{n-1}$ be a genus one normal curve defined over a
field $K$. Suppose we scale the matrix $\Omega$ in Theorem~\ref{thm1}
so that the coefficients of its entries are in $K$. Then $C$ has
Jacobian elliptic curve $y^2 = x^3 - 27 c_4(\Omega) x - 54
c_6(\Omega)$.
\end{Corollary}

\begin{Proof}
  Let $E$ be the elliptic curve $y^2 = x^3 - 27 c_4(\Omega) x - 54
  c_6(\Omega)$.  By Theorem~\ref{thm2} there is an isomorphism $\gamma
  : C \to E$ with $\gamma^*(3 dx/y) = \omega$. Let $\xi_\sigma =
  \sigma(\gamma) \gamma^{-1}$ for $\sigma \in \Gal(\Kbar/K)$. Since
  $3dx/y$ and $\omega$ are both $K$-rational it follows that
  $\xi_\sigma^* (3 dx/y) = 3dx/y$.  This implies, as explained for
  example in \cite[Lemma 2.4]{g1inv}, that $\xi_\sigma : E \to E$
  is a translation map. Then $C$ is the twist
  of $E$ by the class of $\{\xi_\sigma\}$ in $H^1(K,E)$. It follows by
  Theorems 3.6 and 3.8 in \cite[Chapter~X]{Sil} that $C$ is a
  principal homogeneous space under $E$, and $E$ is the Jacobian of
  $C$.
\end{Proof}

\begin{Remark}
  Although we will not need it for the proofs of Theorems~\ref{thm1}
  and~\ref{thm2}, it is natural to ask whether $C \subset \PP^{n-1}$
  is uniquely determined by $\Omega$. The answer is that it is.
  Indeed by the minimal free resolutions in Theorem~\ref{thm1} we can
  recover $\nabla$ from $\Omega$. Then by Euler's identity we obtain
  equations for $\Sec^r C$ where $n-2r = 1$ or $2$.  This then
  determines $\Sec^1 C = C$ by Theorem~\ref{thm:hsec-eqns}(v).
\end{Remark}

\section{Changes of co-ordinates}
\label{sec:change}

We show that the constructions in Section~\ref{sec:stat} behave well
under all changes of co-ordinates.  First we define an action of
$\GL_n$ on the space of all $n \times n$ alternating matrices of
quadratic forms in $x_1, \ldots x_n$. For $g \in \GL_n$ we put
\[ g \star \Omega = g^{-T} \left( \Omega ( \sum_{i=1}^n g_{i1} x_i,
  \ldots, \sum_{i=1}^n g_{in} x_i) \right) g^{-1} \] where $g^{-T}$ is
the inverse transpose of $g$. Since the scalar matrices act trivially,
this could equally be viewed as an action of $\PGL_n$.

\begin{Lemma}
\label{action_on_omega}
Let $C \subset \PP^{n-1}$ and $C' \subset \PP^{n-1}$ be genus one
normal curves. Let $\Omega$ and $\Omega'$ be alternating matrices of
quadratic forms 
that satisfy the conclusions of Theorem~\ref{thm1}, and define
invariant differentials $\omega$ and $\omega'$ on $C$ and $C'$.
If $\gamma : C' \to C$ is an isomorphism given by
\[  (x_1 : \ldots : x_n) \mapsto ( \sum_{i=1}^n g_{i1} x_i :
\ldots : \sum_{i=1}^n g_{in} x_i) \]
for some $g \in \GL_n$ then there exists $\lambda \in K^\times$ 
such that $g \star \Omega = \lambda 
\Omega'$ and $\gamma^* \omega = \lambda^{-1} \omega'$.
\end{Lemma}
\begin{Proof}
  Suppose $n$ is odd, say $n=2r+1$ and $\Sec^r C = \{F = 0\}$. Then
  $\Sec^r C'$ is defined by
  \[ F'(x_1, \ldots, x_n) = F(y_1, \ldots, y_n) \] where $y_j =
  \sum_{i=1}^n g_{ij} x_i$.  By the chain rule
  \[ \nabla(F')(x_1, \ldots,x_n) = \nabla(F) (y_1, \ldots, y_n) \,\,
  g^T. \] Then
  \[ \nabla(F) \Omega = 0 \implies \nabla(F') (g \star \Omega) = 0. \]
  It follows by the uniqueness of minimal free resolutions that $g
  \star \Omega = \lambda \Omega'$ for some $\lambda \in K^\times$.
  The case $n$ is even is similar.

  We also have $\gamma^* \omega = \mu \omega'$ for some $\mu \in
  K^\times$.  If $y_j = \sum_{i=1}^n g_{ij} x_i$ then
  \[ y_s^2 d(y_r/y_s) = \sum_{i,j=1}^n g_{ir} g_{js} x_j^2
  d(x_i/x_j). \] Dividing by $\gamma^* \omega = \mu \omega'$ gives
  \[ \Omega(y_1, \ldots, y_n) = \mu^{-1} g^T \Omega'(x_1,\ldots,x_n)
  g. \] Hence $g \star \Omega = \mu^{-1} \Omega'$ and so $\mu =
  \lambda^{-1}$.
\end{Proof}

\begin{Lemma}
\label{leminv}
The polynomials $c_4$ and $c_6$ are invariants for the action of
$\GL_n$, i.e. $c_4( g \star \Omega) = c_4(\Omega)$ and $c_6( g \star
\Omega) = c_6(\Omega)$ for all $g \in \GL_n$.
\end{Lemma}

\begin{Proof}
  Let $\Omega' = g \star \Omega$, i.e.
  \[ \Omega'_{ij}(x_1, \ldots, x_n) = \sum_{a,b=1}^n (g^{-1})_{ai}
  (g^{-1})_{bj} \Omega_{ab}(y_1, \ldots, y_n) \] where $y_j =
  \sum_{i=1}^n g_{ij} x_i$. Direct calculation using~\eqref{MN} shows
  that
\begin{align*}
  M'_{ij}(x_1, \ldots, x_n) &= \sum_{a,b=1}^n (g^{-1})_{ai}
  (g^{-1})_{bj} M_{ab}(y_1, \ldots, y_n), \\
  N'_{ijk}(x_1, \ldots, x_n) &= \sum_{a,b,c=1}^n (g^{-1})_{ai}
  (g^{-1})_{bj} (g^{-1})_{ck} N_{abc}(y_1, \ldots, y_n).
\end{align*}
Then
\begin{align*}
  \frac{\partial^2 M'_{ij}}{\partial x_r \partial x_s} &=
  \sum_{a,b,c,d=1}^n (g^{-1})_{ai} (g^{-1})_{bj} g_{rc} g_{sd}
  \frac{\partial^2 M_{ab}}{\partial x_c \partial x_d}, \\
  \frac{\partial^2 M'_{rs}}{\partial x_i \partial x_j} &=
  \sum_{A,B,C,D=1}^n (g^{-1})_{Cr} (g^{-1})_{Ds} g_{iA} g_{jB}
  \frac{\partial^2 M_{CD}}{\partial x_A \partial x_B}.
\end{align*}
Multiplying these together and summing gives
\[ \sum_{i,j,r,s=1}^n \frac{\partial^2 M'_{ij}}{\partial x_r \partial
  x_s} \frac{\partial^2 M'_{rs}}{\partial x_i \partial x_j} =
\sum_{a,b,c,d=1}^n \frac{\partial^2 M_{ab}}{\partial x_c \partial x_d}
\frac{\partial^2 M_{cd}}{\partial x_a \partial x_b}.  \] Thus
$c_4(\Omega') = c_4(\Omega)$. A similar argument shows that
$c_6(\Omega') = c_6(\Omega)$.
\end{Proof}

The following corollary shows that to prove Theorems~\ref{thm1}
and~\ref{thm2} for a fixed value of $n$, it suffices to prove them for
a family of curves covering the $j$-line.

\begin{Corollary}
\label{howc4c6change}
Let $\Omega_1$ and $\Omega_2$ correspond to pairs $(C_1,\omega_1)$ and
$(C_2,\omega_2)$. If there is an isomorphism $\gamma : C_1 \to C_2$
with $\gamma^* \omega_2 = \lambda \omega_1$ then $c_4(\Omega_1) =
\lambda^4 c_4(\Omega_2)$ and $c_6(\Omega_1) = \lambda^6
c_6(\Omega_2)$.
\end{Corollary}
\begin{Proof} After composing the isomorphism $\gamma$ with a
  translation map, we may suppose it is given by a change of
  co-ordinate on $\PP^{n-1}$.  The case $\lambda = 1$ is immediate
  from Lemmas~\ref{action_on_omega} and~\ref{leminv}. In general we
  use that $c_4$ and $c_6$ are homogeneous polynomials 
  of degrees $4$ and $6$.
\end{Proof}

\section{Curves of small degree}

We compare our general formula for the Jacobian with the formulae
previously known for genus one normal curves of degrees $3$, $4$ and
$5$.

For curves of degrees $3$ and $4$ it is easy to write down a matrix
$\Omega$ satisfying the conclusions of Theorems~\ref{thm1}
and~\ref{thm2}(i).  Indeed for $C = \{F(x_1,x_2,x_3) = 0\} \subset
\PP^2$ a plane cubic we put
\[ \Omega = \begin{pmatrix} \vspace{0.5ex} 0 & \frac{\partial F
  }{\partial x_3} & -\frac{\partial F }{\partial x_2} \\
  \vspace{0.5ex}
  -\frac{\partial F }{\partial x_3} & 0 & \frac{\partial F }{\partial x_1} \\
  \frac{\partial F }{\partial x_2} & -\frac{\partial F }{\partial x_1}
  & 0 \\ \end{pmatrix}, \] and for $C = \{F_1 = F_2 = 0\} \subset
\PP^3$ a quadric intersection we let $\Omega$ be the $4 \times 4$
alternating matrix with entries
\[ \Omega_{ij} = \frac{\partial F_1 }{\partial x_k} \frac{\partial F_2
}{\partial x_l} - \frac{\partial F_1 }{\partial x_l} \frac{\partial
  F_2 }{\partial x_k}, \] where $(i,j,k,l)$ is an even permutation of
$(1,2,3,4)$.  To prove Theorem~\ref{thm2}(ii) in these cases we may
check by direct computation that $c_4(\Omega)$ and $c_6(\Omega)$ are
the classical invariants of a ternary cubic or quadric intersection,
as scaled in \cite[Section 7]{g1inv}.  We note that these are
polynomials of degrees 4 and 6 in the coefficients of $F$,
respectively of degrees 8 and 12 in the coefficients of
$F_1$ and $F_2$.

As described for example in \cite[Section 4]{5desc}, a genus one
normal curve of degree $n=5$ is defined by the $4 \times 4$ Pfaffians
$p_1, \ldots, p_5$ of a $5 \times 5$ alternating matrix of linear
forms on $\PP^4$. We call the matrix of linear forms $\Phi$ a {\em
  genus one model} of degree $5$, and note that there is a natural
action of $\GL_5 \times \GL_5$ on the space of all such models. It is
shown in \cite[Proposition~VIII.2.5]{Hulek} that the secant variety
$\Sec^2 C$ is a hypersurface of degree $5$ with equation $F=0$ where
$F$ is the determinant of the Jacobian matrix of $p_1, \ldots,
p_5$. In \cite[Section 7]{invenqI} we proved that there is a degree
$5$ covariant $\Omega$ satisfying the conclusions of
Theorems~\ref{thm1} and~\ref{thm2}(i). We gave an explicit formula for
this covariant in~\cite[Section 2]{sqrfree}.

We claim that $c_4(\Omega)$ and $c_6(\Omega)$ are invariants for the
action of $\SL_5 \times \SL_5$. For the action of $\SL_5$ via changes
of co-ordinates on $\PP^4$ this follows from Lemma~\ref{leminv}. For
the action of $\SL_5$ via $\Phi \mapsto A \Phi A^T$ it turns out that
the coefficients of the entries of $\Omega$ are already
invariants. Since $\Omega$ is a covariant of degree $5$, the
invariants $c_4(\Omega)$ and $c_6(\Omega)$ have degrees $20$ and $30$
in the coefficients of the entries of $\Phi$.  Computing a single
numerical example (to check the scaling) shows that $c_4(\Omega)$ and
$c_6(\Omega)$ are the same as the invariants $c_4(\Phi)$ and
$c_6(\Phi)$ constructed in \cite{g1inv}.

\section{Minimal free resolutions}
\label{sec:outline}

Let $C \subset \PP^{n-1}$ be a genus one normal curve of degree $n \ge
3$.  Let $\Omega$ be an $n \times n$ alternating matrix of quadratic
forms in $x_1, \ldots, x_n$.  In Sections~\ref{sec:pf1}
and~\ref{sec:checkhyp} we exhibit $\Omega$ satisfying the following
three hypotheses.

\renewcommand{\theenumi}{{\bf H\arabic{enumi}}}

\medskip

\begin{enumerate}
\item If $n-2r \ge 1$ and $f \in I(\Sec^r C)$ then $\sum_{i=1}^n
  \frac{\partial f}{\partial x_i} \Omega_{ij} \in I(\Sec^r C)$ for all
  $j$.
\item \smallskip If $n - 2r = 2$ and $\Sec^r C = \{F_1 = F_2 = 0\}$
  then $\sum_{i,j=1}^n \frac{\partial F_1}{\partial x_i} \Omega_{ij}
  \frac{\partial F_2}{\partial x_j} = 0$.
\item \smallskip If $n-2r \ge 1$ then there exists $P \in \Sec^r C$
  with $\rank \Omega(P) = 2r$.
\end{enumerate}

\medskip

\renewcommand{\theenumi}{\roman{enumi}}

In this section we prove:

\begin{Theorem}
\label{thm3}
Let $\Omega$ be an $n \times n$ alternating matrix of quadratic forms,
satisfying the hypotheses {\bf (H1)}, {\bf (H2)} and {\bf (H3)}. Then
there is a minimal free resolution as described in Theorem~\ref{thm1}.
\end{Theorem}

The next two propositions are proved in Section~\ref{sec:hsec}.  By
abuse of notation we write $P$ both for a point in $\PP^{n-1}$ and for 
a vector of length $n$ representing this point.

\begin{Proposition}
\label{lem:tangent}
If $n-2r \ge 1$ and $P = \sum_{i=1}^r \xi_i P_i$ for some $P_1,
\ldots, P_r \in C$ distinct and $\xi_1, \ldots, \xi_r \not=0$ then the
tangent space $T_P \Sec^r C$ is the linear span of the tangent lines
$T_{P_1} C, \ldots, T_{P_r} C$.
\end{Proposition}

\pagebreak

\begin{Proposition}
\label{lem:codim3}
Let $\nabla(F)$ and $\nabla(F_1,F_2)$ be as defined in
Theorem~\ref{thm1}.
\begin{enumerate}
\item If $n-2r=1$ and $\Sec^r C = \{ F = 0 \}$ then the entries of
  $\nabla(F)$ define a variety in $\PP^{n-1}$ of codimension $3$.
\item If $n-2r=2$ and $\Sec^r C = \{ F_1 = F_2 = 0\}$ then the $2
  \times 2$ minors of $\nabla(F_1,F_2)$ define a variety in
  $\PP^{n-1}$ of codimension $3$.
\end{enumerate}
\end{Proposition}
\begin{Proof}
  (i) Theorem~\ref{thm:hsec-eqns} tells us that $\Sec^r C$ has
  singular locus $\Sec^{r-1} C$, and that this has codimension $3$. \\
  (ii) This is proved in Section~\ref{sec:codim3}.
\end{Proof}

We start the proof of Theorem~\ref{thm3} with the following lemma.

\begin{Lemma}
\label{lem2}
Let $C \subset \PP^{n-1}$ be a genus one normal curve. Suppose that
$n-2r \ge 1$ and $\ell_1, \ldots, \ell_n$ are linear forms in $x_1,
\ldots, x_n$ such that
\begin{equation}
\label{linsyz}
\sum_{i=1}^n \ell_i \frac{\partial f}{\partial x_i} \in I(\Sec^r C) 
\qquad \text{ for all } f \in I(\Sec^r C). 
\end{equation}
Then there exists $\lambda \in K$ such that $\ell_i = \lambda x_i$ for
all $1 \le i \le n$.
\end{Lemma}

\begin{Proof}
  The coefficients of $\ell_1, \ldots, \ell_n$ form an $n \times n$
  matrix.  Let $V \subset \Mat_n(K)$ be the subspace of all solutions
  to~\eqref{linsyz}.  We must show that $V$ consists only of scalar
  matrices.  Let $E$ be the Jacobian of $C$. Translation by $T \in
  E[n]$ is an automorphism of $C$ that extends to an automorphism of
  $\PP^{n-1}$, say given by a matrix $M_T$. Now $V$ is stable under
  conjugation by each $M_T$.  By considering the standard
  representation of the Heisenberg group (see for example
  \cite[Section 3]{pfpres}) it follows that $V$ has a basis $\{ M_T :
  T \in X \}$ for some subset $X \subset E[n]$.

  We suppose for a contradiction that $M_T \in V$ for some $0_E \not=
  T \in E[n]$. Then translation by $T$ on $C$ extends to an
  automorphism of $\PP^{n-1}$ that sends each point $P \in \Sec^r C$
  to a point in the tangent space $T_P \Sec^r C$.  Let $H$ be the
  divisor of a hyperplane section on $C$. For $D$ an effective divisor
  on $C$ we write $\overline{D} \subset \PP^{n-1}$ for the linear
  subspace cut out by $\LL(H-D) \subset \LL(H)$. For example, if $D$
  is a sum of distinct points on $C$ then $\overline{D}$ is the linear
  span of these points.  We also write $D_T$ for $D$ translated by
  $T$. We choose $D = P_1 + \ldots + P_r$ an effective divisor of
  degree $r$ such that
  \begin{enumerate}
  \item $P_1, \ldots, P_r \in C$ are distinct,
  \item $D$ and $D_T$ have disjoint support,
  \item $2 D + D_T \not\sim H$.
  \end{enumerate}
  Proposition~\ref{lem:tangent} shows that for generic $P \in
  \overline{D}$ we have $T_P \Sec^r C = \overline{2D}$. It follows
  from our assumption $M_T \in V$ that $\overline{D_T} \subset
  \overline{2D}$, equivalently $\LL(H - 2 D) \subset \LL(H - D_T)$.
  Then by (ii) we have
  \[ \LL(H - 2D) = \LL(H - 2 D) \cap \LL(H - D_T) = \LL(H - 2D -
  D_T). \] However by (iii) and the Riemann-Roch theorem these spaces
  do not have the same dimension.  Indeed, since $r \ge 1$ and $n - 2r
  \ge 1$ we have
  \[ \dim \LL(H - 2D) = n - 2 r \not= \max(n - 3 r,0) = \dim \LL(H -
  2D-D_T). \] This is the required contradiction.
\end{Proof}

We show that the resolution in Theorem~\ref{thm1} is a complex.

\begin{Lemma}
\label{lem:complex}
Let $C \subset \PP^{n-1}$ be a genus one normal curve, 
and let $\Omega$ be an alternating matrix of quadratic forms
satisfying the hypotheses {\bf (H1)} and {\bf (H2)}.
\begin{enumerate}
\item If $n=2r+1$ and $\Sec^r C = \{ F = 0 \}$ then
\[ \sum_{i=1}^n \frac{\partial F}{\partial x_i} \Omega_{ij} = 0 \qquad
\text{ for all } 1 \le j \le n. \]
\item  If $n=2r+2$ and $\Sec^r C = \{ F_1 = F_2 = 0 \}$ then
\[ \sum_{i=1}^n \frac{\partial F_1}{\partial x_i} \Omega_{ij} = 
   \sum_{i=1}^n \frac{\partial F_2}{\partial x_i} \Omega_{ij} = 0 \qquad
\text{ for all } 1 \le j \le n. \]
\end{enumerate}
\end{Lemma}
\begin{Proof}
  (i) By the hypothesis {\bf (H1)} we have
\[  \sum_{i=1}^n \frac{\partial F}{\partial x_i} \Omega_{ij} = \ell_j F 
\qquad \text{ for all } 1 \le j \le n, \] 
for some linear forms $\ell_1, \ldots, \ell_n$. We multiply by 
$\frac{\partial F}{\partial x_j}$ and sum over $j$. Since $\Omega$ is 
alternating the left hand side is zero. Therefore
\begin{equation*}
\sum_{j=1}^n \ell_j \frac{\partial F}{\partial x_j} = 0. 
\end{equation*}
By Lemma~\ref{lem2} and Euler's identity it follows that $\ell_1 =
\ldots = \ell_n = 0$ as required. \\
(ii) By the hypothesis {\bf (H1)} we have
\begin{equation}
\label{eqn:new}
  \sum_{i=1}^n \frac{\partial F_1}{\partial x_i} \Omega_{ij} = \ell_j F_1
 + m_j F_2 \qquad \text{ for all } 1 \le j \le n, 
\end{equation}
for some linear forms $\ell_1, \ldots, \ell_n$ and $m_1,\ldots,m_n$.
We multiply by $\frac{\partial F_1}{\partial x_j}$ and sum over~$j$.
Since $\Omega$ is alternating the left hand side is zero.  Since $F_1$
and $F_2$ are forms defining a variety of codimension $2$ they must be
coprime. Therefore
\begin{equation*}
  \sum_{j=1}^n \ell_j \frac{\partial F_1}{\partial x_j} = \xi F_2 
  \quad \text{ and }
  \quad \sum_{j=1}^n m_j \frac{\partial F_1}{\partial x_j} = -\xi F_1
\end{equation*}
for some $\xi \in K$.  If instead we multiply~\eqref{eqn:new} by
$\frac{\partial F_2}{\partial x_j}$ and sum over $j$ then using the
hypothesis {\bf (H2)} we find that
\begin{equation*}
  \sum_{j=1}^n \ell_j \frac{\partial F_2}{\partial x_j} = \eta F_2 \quad 
  \text{ and } \quad \sum_{j=1}^n m_j \frac{\partial F_2}{\partial x_j}
  = -\eta F_1
\end{equation*}
for some $\eta \in K$.

By Lemma~\ref{lem2} there exist $\lambda, \mu \in K$ such that $\ell_i
= \lambda x_i$ and $m_i = \mu x_i$ for all $1 \le i \le n$.  By
Euler's identity and the linear independence of $F_1$ and $F_2$ it
follows that $\lambda = \mu = 0$. Therefore
\[ \sum_{i=1}^n \frac{\partial F_1}{\partial x_i} \Omega_{ij} = 0
\qquad \text{ for all } 1 \le j \le n. \] 
The corresponding result for $F_2$ follows by symmetry.
\end{Proof}

To complete the proof of Theorem~\ref{thm3} we must 
show that the complex is exact.  First we need some
linear algebra.  If $B$ is an $n \times n$ matrix and $S \subset \{1,
\ldots, n\}$ then we write $B^S$ for the $(n-|S|) \times (n-|S|)$
matrix obtained by deleting the rows and columns indexed by $S$. The
Pfaffian $\pf(M)$ of an alternating matrix $M$ is a polynomial in the
matrix entries with the property that $\det(M) = \pf(M)^2$.

\begin{Lemma}
\label{la}
\begin{enumerate}
\item Let $A$ be a $1 \times n$ matrix and $B$ an $n \times n$
  alternating matrix over a field $K$. Suppose that $\rank A = 1$,
  $\rank B = n-1$ and $AB = 0$.  Then there exists $\lambda \in
  K^\times$ such that
\[ (-1)^i \pf (B^{\{i\}}) = \la a_i \]
for all $1 \le i \le n$. 
\item Let $A$ be a $2 \times n$ matrix and $B$ an $n \times n$
  alternating matrix over a field $K$. Suppose that $\rank A = 2$,
  $\rank B = n-2$ and $AB = 0$.  Then there exists $\lambda \in
  K^\times$ such that
  \[ (-1)^{i+j} \pf (B^{\{i,j\}}) = \la (a_{1i} a_{2j} - a_{1j}
  a_{2i}) \] for all $1 \le i < j \le n$.
\end{enumerate}
\end{Lemma}
\begin{Proof}
  (i) It is well known that the vector with $i$th entry $(-1)^i \pf
  (B^{\{i\}})$ belongs to the kernel of $B$. See for example
  \cite[Section 3.4]{CMrings}.  Since $\rank B = n-1$ this vector
  is non-zero and the kernel is $1$-dimensional. The result follows. \\
  (ii) We first claim there exist $\lambda_1, \ldots, \lambda_n \in K$
  such that
\[ (-1)^{i+j} \pf (B^{\{i,j\}}) = \left\{
  \begin{array}{ll} \la_i (a_{1i} a_{2j} - a_{1j} a_{2i}) & \text{ if
    } i<j, \\ -\la_i (a_{1i} a_{2j} - a_{1j} a_{2i}) & \text{ if }
    i>j. \end{array} \right. \] 
Indeed taking $a_{2i}$ times the first row of $A$ minus $a_{1i}$ times
the second row of $A$, gives a non-zero vector in the kernel of
$B^{\{i\}}$.  If $\rank B^{\{i\}} = n-2$ then we argue as in (i).
Otherwise we can simply take $\lambda_i = 0$. This proves the claim.

Now let $C = (a_{1i} a_{2j} - a_{1j} a_{2i})_{i,j=1,\ldots,n}$ and let
$D$ be the diagonal matrix with entries $\lambda_1, \ldots,
\lambda_n$.  We must show that if $CD = DC$ then $CD$ is a scalar
multiple of $C$. More generally this is true for any rank $2$
alternating matrix $C$ and diagonal matrix $D$. Indeed we may re-order
the rows and columns so that the diagonal entries of $D$ which are
equal are grouped together.  Then $C$ is in block diagonal form. Since
$C$ is alternating of rank $2$, exactly one of these blocks is
non-zero. The result is then clear.
\end{Proof}

\begin{Lemma}
\label{lem:pfs}
Let $C \subset \PP^{n-1}$ be a genus one normal curve,
and let $\Omega$ be an alternating matrix of quadratic forms
satisfying the hypotheses {\bf (H1)}, {\bf (H2)} and {\bf (H3)}.
\begin{enumerate}
\item If $n=2r+1$ and $\Sec^r C = \{F = 0\}$ then the $(n-1) \times
  (n-1)$ Pfaffians of $\Omega$ are (scalar multiples of) the partial
  derivatives of $F$.
\item If $n=2r+2$ and $\Sec^r C = \{ F_1 = F_2 = 0 \}$ then the $(n-2)
  \times (n-2)$ Pfaffians of $\Omega$ are (scalar multiples of) the $2
  \times 2$ minors of $\nabla(F_1,F_2)$.
\end{enumerate}
\end{Lemma}
\begin{Proof}
  We apply Lemma~\ref{la} over the function field
  $K(x_1, \ldots,x_n)$. \\
  (i) By Lemma~\ref{lem:complex} we have $\sum_{i=1}^n \frac{\partial
    F}{\partial x_i} \Omega_{ij} = 0$. By the hypothesis {\bf (H3)}
  the generic rank of $\Omega$ is $n-1$. So by Lemma~\ref{la}(i) there
  exists $\la \in K(x_1, \ldots, x_n)$ such that
  \[ (-1)^i \pf (\Omega^{\{i\}}) = \la \frac{\partial F}{\partial x_i}
  \qquad \text{ for all $1 \le i \le n$. } \] Since $\pf
  (\Omega^{\{i\}})$ and $\tfrac{\partial F}{\partial x_i}$ are forms
  of degree $n-1$, we can write $\lambda = u/v$ where $u$ and $v$ are
  coprime forms of the same degree. Then $v$ divides $\tfrac{\partial
    F}{\partial x_i}$ for all $i$, and so must be a constant by
  Proposition~\ref{lem:codim3}(i). Therefore $\lambda$ is
  a constant. \\
  (ii) By Lemma~\ref{lem:complex} we have $\sum_{i=1}^n \frac{\partial
    F_1}{\partial x_i} \Omega_{ij} = \sum_{i=1}^n \frac{\partial
    F_2}{\partial x_i} \Omega_{ij} = 0$. By the hypothesis {\bf (H3)}
  the generic rank of $\Omega$ is $n-2$. So by Lemma~\ref{la}(ii)
  there exists $\la \in K(x_1, \ldots, x_n)$ such that
  \[ (-1)^{i+j} \pf (\Omega^{\{i,j\}}) = \la
  \frac{\partial(F_1,F_2)}{\partial(x_i,x_j)} \qquad \text{ for all $1
    \le i<j \le n$.} \] Since $\pf (\Omega^{\{i,j\}})$ and
  $\tfrac{\partial(F_1,F_2)}{\partial(x_i,x_j)}$ are forms of degree
  $n-2$, we can write $\lambda = u/v$ where $u$ and $v$ are coprime
  forms of the same degree. Then $v$ divides
  $\frac{\partial(F_1,F_2)}{\partial(x_i,x_j)}$ for all $i, j$, and so
  must be a constant by Proposition~\ref{lem:codim3}(ii). Therefore
  $\lambda$ is a constant.
\end{Proof}

Let $R = K[x_1, \ldots, x_n]$. Consider a complex of graded free
$R$-modules
\begin{equation}
\label{complex}
0 \to F_m \stackrel{\varphi_m}{\ra} F_{m-1} \ra \ldots 
\ra F_1 \stackrel{\varphi_1}{\ra} F_0. 
\end{equation}
We write $V_k \subset \PP^{n-1}$ for the subvariety defined by the
$r_k \times r_k$ minors of $\varphi_k$ where $r_k = \rank
(\varphi_k)$.  The Buchsbaum-Eisenbud acyclicity criterion (see
\cite[Theorem~1.4.13]{CMrings} or \cite[Theorem~20.9]{Eisenbud})
states that~\eqref{complex} is exact if and only if $\rank F_k = \rank
\varphi_k + \rank \varphi_{k+1}$ and $\codim V_k \ge k$ for all $1 \le
k \le m$.

\medskip

\begin{ProofOf}{Theorem~\ref{thm3}}
  We already saw in Lemma~\ref{lem:complex} that the resolution in
  Theorem~\ref{thm1} is a complex.  We must prove it is exact. If $n$
  is odd then the free $R$-modules have ranks $1,n,n,1$ and the maps
  have ranks $1,n-1,1$.  If $n$ is even then the free $R$-modules have
  ranks $2,n,n,2$ and the maps have ranks $2,n-2,2$.  By
  Lemma~\ref{lem:pfs} we have $V_1 = V_2 = V_3$ and
  Proposition~\ref{lem:codim3} shows that this variety has codimension
  $3$. We now apply the Buchsbaum-Eisenbud acyclicity criterion.
\end{ProofOf}

\section{A basis-free construction}
\label{sec:pf1}

The results of Section~\ref{sec:change}
show that for the proof of
Theorems~\ref{thm1} and~\ref{thm2} we are free to make changes of
co-ordinates on $\PP^{n-1}$.  Since we are working over an
algebraically closed field we can therefore reduce to the following
situation.  Let $E$ be the elliptic curve
\begin{equation*}
y^2 + a_1 xy + a_3 y = x^3 + a_2 x^2 + a_4 x + a_6 
\end{equation*}
with point at infinity $0_E$ and invariant differential
\[\omega = dx/(2y + a_1 x + a_3) = dy/(3 x^2 + 2a_2 x+a_4-a_1y).\] 
Let $C \subset \PP^{n-1}$ be the image of $E$ embedded via the
complete linear system $|n.0_E|$. 
The embedding depends on a choice of basis for
the Riemann-Roch space $\LL(n.0_E)$, but the only effect
of changing this is to make a change of co-ordinates on $\PP^{n-1}$.
In this section we define a linear map $\Omega :
\wedge^2 \LL(n.0_E) \to S^2 \LL(n.0_E)$. In the next section we show
that the corresponding alternating matrix of quadratic forms satisfies
the hypotheses {\bf (H1)}, {\bf (H2)} and {\bf (H3)}.
 
For $f \in \LL(n.0_E)$ we put $\dot{f} = df/\omega \in
\LL((n+1).0_E)$.  Motivated by~\eqref{firstomega} we define a linear
map
\begin{align*}
  A : \wedge^2 \LL(n.0_E) & \to S^2 \LL((n+1).0_E) \\
  f \wedge g & \mapsto f \otimes \dot{g} - g \otimes \dot{f}.
\end{align*}
\begin{Lemma} 
\label{lem:1}
Let $f,g \in \LL(n.0_E)$. Then the rational function on $E \times E$
given by
\[ (P,Q) \mapsto \frac{y_P + y_Q + a_1 x_Q + a_3}{x_P - x_Q} (f(Q)
g(P) - f(P) g(Q)) \] belongs to $\LL((n+1).0_E) \otimes
\LL((n+1).0_E)$.
\end{Lemma}

\begin{Proof}
(i) If we fix $Q=(x_Q,y_Q)$ then as rational functions of $P = (x,y)$,
\[ \frac{y + y_Q + a_1 x_Q + a_3}{x - x_Q} \in \LL(0_E + Q) \quad
\text{ and } \quad f(Q) g - g(Q) f \in \LL(n.0_E - Q). \]
Therefore the product belongs to $\LL((n+1).0_E)$. \\
(ii) If we fix $P=(x_P,y_P)$ then as rational functions of $Q =
(x,y)$,
\[ \frac{y_P + y + a_1 x + a_3}{x_P - x} \in \LL(0_E + P) \quad \text{
  and } \quad g(P) f - f(P) g \in \LL(n.0_E - P). \] Therefore the
product belongs to $\LL((n+1).0_E)$.
\end{Proof}

We define a second linear map
\begin{align*}
  B : \wedge^2 \LL(n.0_E) & \to S^2 \LL((n+1).0_E) \\
  f \wedge g & \mapsto \frac{y_P + y_Q + a_1 x_Q + a_3}{x_P - x_Q}
  (f(Q) g(P) - f(P) g(Q)) \bigg|_{P=Q}
\end{align*}
where $|_{P=Q}$ is our notation for the natural map 
\[\LL((n+1).0_E) \otimes \LL((n+1).0_E) \to S^2 \LL((n+1).0_E).\]

We show that $A$ and $B$ both represent the invariant differential
$\omega$, in the sense of Theorem~\ref{thm2}(i).

\begin{Lemma} 
\label{lem:omega}
As rational functions on $E$ we have
\[  A(f \wedge g) = B(f \wedge g) = f \dot{g} - g \dot{f} = 
\frac{f dg - g df }{\omega}. \]
\end{Lemma}
\begin{Proof}
This is clear for $A$. For $B$ we apply l'H\^opital's rule to get
\[  \frac{ f(Q) g - g(Q) f }{x - x_Q} \bigg|_{P=Q} =  
\frac{ f(Q) \dot{g} - g(Q) \dot{f} }{ \dot{x} } \bigg|_{P=Q}, \]
and then use that $\dot{x} = 2 y + a_1 x + a_3$.
\end{Proof}

If we pick bases for $\LL(n.0_E)$ and $\LL((n+1).0_E)$ then $A$ and
$B$ are (represented by) $n \times n$ alternating matrices of
quadratic forms in $n+1$ variables. However the matrix $\Omega$ we
seek is an $n \times n$ alternating matrix of quadratic forms in $n$
variables. It turns out that the correct choice of $\Omega$ is a
linear combination of $A$ and $B$.

We may expand rational functions on $E$ as Laurent power series in the
local parameter $t=x/y$ at $0_E$. Let $\phi$ be the linear map that
reads off the coefficient of $t^{-n-1}$. There are exact sequences
\[ 0 \to \LL(n.0_E) \to \LL((n+1).0_E) \stackrel{\phi}{\to} K \to 0 \]
and
\begin{equation}
  \label{lastcol}
  0 \to S^2 \LL(n.0_E) \to S^2 \LL((n+1).0_E) \stackrel{\phi_2}{\to}  
  \LL((n+1).0_E) \to 0 
\end{equation}
where $\phi_2(f \otimes g) = \phi(f) g + \phi(g) f$.

\begin{Lemma} 
\label{lem:3}
Let $f,g \in \LL(n.0_E)$ be rational functions whose coefficients of
$t^{-n}$ (when expanded as Laurent power series in $t$) are $0, 1$
respectively. Then
\[  \phi_2(A(f \wedge g)) = nf \quad \text{ and } 
\quad \phi_2(B(f \wedge g)) = 2f. \]
\end{Lemma}

\begin{Proof}
  (i) We have $x = t^{-2} + \ldots$ and $y = t^{-3} + \ldots$. Then
  $\dot{x} = 2y + a_1 x + a_3 = 2 t^{-3} + \ldots$ and $\dot{y} = 3
  x^2 + 2 a_2 x + a_4 - a_1 y = 3 t^{-4} + \ldots$. Writing $g$ as a
  polynomial in $x$ and $y$ it follows that $g = t^{-n} + \ldots$
  and $\dot{g} = n t^{-n-1} + \ldots$. Therefore $\phi(f) = \phi(g) =
  \phi(\dot{f})=0$ and $\phi( \dot{g} ) = n$. We compute
  \[ \phi_2(A(f \wedge g)) = \phi_2( f \otimes \dot{g} - g \otimes
  \dot{f} ) = nf. \] (ii) If we fix $Q=(x_Q,y_Q)$ then as rational
  functions of $P = (x,y)$,
  \[ \frac{y + y_Q + a_1 x_Q + a_3}{x - x_Q} = t^{-1} + \ldots \quad
  \text{ and } \quad f(Q) g - g(Q) f = f(Q) t^{-n} + \ldots \] with
  product $f(Q) t^{-n-1} + \ldots$.

  If we fix $P=(x_P,y_P)$ then as rational functions of $Q = (x,y)$,
  \[ \frac{y_P + y + a_1 x + a_3}{x_P - x} = -t^{-1} + \ldots \quad
  \text{ and } \quad g(P) f - f(P) g = -f(P) t^{-n} + \ldots \] with
  product $f(P) t^{-n-1} + \ldots$.  In both cases the leading
  coefficient is $f$. Adding these together gives $\phi_2(B(f \wedge
  g)) = 2f$.
\end{Proof}

\begin{Corollary}
\label{cor:4}
Let $\Omega = n B - 2 A$. Then $\Omega$ is a linear map $\wedge^2
\LL(n.0_E) \to S^2 \LL(n.0_E)$.
\end{Corollary}
\begin{Proof}
This follows from Lemma~\ref{lem:3} and the exact sequence~\eqref{lastcol}.
\end{Proof}

\section{Proof of Theorem~\ref{thm1}}
\label{sec:checkhyp}

If we pick a basis for $\LL(n.0_E)$ then the linear map defined in
Corollary~\ref{cor:4} is represented by an $n \times n$ alternating
matrix of quadratic forms in $n$ variables. In this section we
complete the proof of Theorem~\ref{thm1} by showing that this matrix
$\Omega$ satisfies the hypotheses {\bf (H1)}, {\bf (H2)} and {\bf
  (H3)}, as stated at the start of Section~\ref{sec:outline}.

For $0_E \not= P \in E$ we write $\vP$ and $d\vP$ for the linear maps
$f \mapsto f(P)$ and $f \mapsto \dot{f}(P)$ in the dual space
$\LL(n.0_E)^*$.  For example, if $\LL(n.0_E)$ has basis $1,x,y,x^2,xy,
\ldots$ then
\begin{align*}
  \vP &= (1,x_P,y_P,x_P^2,x_P y_P, \ldots ), \\
  d \vP & = (0,2y_P + a_1 x_P + a_3,3 x_P^2 + 2 a_2 x_P + a_4 - a_1
  y_P, \ldots ).
\end{align*}
We note that $[\vP]$ is a point on $C \subset \PP^{n-1} =
\PP(\LL(n.0_E)^*)$, with tangent line passing through $[d\vP]$. The
square brackets indicate that we are taking the $1$-dimensional
subspaces spanned by these vectors, i.e. the corresponding points in
projective space.  For $0_E \not= Q \in E$ we likewise define $\vQ$
and $d\vQ$.

For $P,Q \in E$ let $\lambda_{P,Q}$ be the slope of the chord (or
tangent line if $P=Q$) joining $P$ and $Q$. In the following lemma the
vectors $\vP, \vQ, d\vP, d\vQ$ in $\LL(n.0_E)^*$ are extended to
$\LL((n+1).0_E)^*$ using exactly the same definition. Evaluating $A$
or $B$ at a linear combination $\xi \vP + \eta \vQ$ gives an element
of $( \wedge^2 \LL(n.0_E))^* = \wedge^2( \LL(n.0_E)^* )$.

\begin{Lemma}
\label{lem:5}
Let $0_E \not= P,Q \in E$ and $\xi,\eta \in K$. Then
\begin{align*}
  A(\xi \vP + \eta \vQ) &= \xi^2 (\vP \wedge d\vP) + \xi \eta ( \vP
  \wedge d \vQ
  + \vQ \wedge d\vP ) + \eta^2 ( \vQ \wedge d\vQ ), \\
  B(\xi \vP + \eta \vQ) &= \xi^2 (\vP \wedge d\vP) + \xi \eta (
  \lambda_{Q,-P} - \lambda_{P,-Q} ) (\vP \wedge \vQ) + \eta^2 ( \vQ
  \wedge d\vQ ).
\end{align*}
\end{Lemma}
\begin{Proof}
  (i) For $f,g \in \LL(n.0_E)$ we compute
  \[ A(\vP)(f \wedge g) = (f \dot{g} - g \dot{f})(P) = (\vP \wedge
  d\vP) (f \wedge g). \]
  The formula for $A(\xi \vP + \eta \vQ)$ follows by bilinearity. \\
  (ii) For $f,g \in \LL(n.0_E)$ we write
  \[ B(\xi \vP + \eta \vQ)(f \wedge g) = \xi^2 B_0 + \xi \eta B_1 +
  \eta^2 B_2 \] By Lemma~\ref{lem:omega} we have
  \begin{align*}
    B_0 = (f \dot{g} - g \dot{f})(P) &= (\vP \wedge d\vP) (f \wedge g), \\
    B_2 = (f \dot{g} - g \dot{f})(Q) &= (\vQ \wedge d\vQ) (f \wedge
    g).
\end{align*}
Since for $s,t \in \LL((n+1).0_E)$ we have
\begin{align*}
  (s \otimes t)(\xi \vP + \eta \vQ) &=
  s(\xi \vP + \eta \vQ) t(\xi \vP + \eta \vQ) \\
  &= \xi^2 s(P) t(P) + \xi \eta (s(P) t(Q) + s(Q) t(P)) + \eta^2 s(Q)
  t(Q),
\end{align*}
it follows from the definition of $B$ that
\begin{align*}
  B_1 &= \lambda_{P,-Q} (f(Q) g(P) - f(P) g(Q)) + \lambda_{Q,-P}
  (f(P) g(Q) - f(Q) g(P)) \\
  &= (\lambda_{Q,-P} - \lambda_{P,-Q}) (\vP \wedge \vQ)(f \wedge g).
\end{align*}
\end{Proof}

We pick a basis for $\LL(n.0_E)$, so that now $\Omega(\vP)$ is an $n
\times n$ alternating matrix, and $\vP$, $\vQ$, $d\vP$, $d\vQ$ are
column vectors.

\begin{Lemma}
\label{lem:6}
Let $0_E \not= P_1, \ldots, P_r \in E$ distinct and $\xi_1, \ldots,
\xi_r \in K$.  Then
\[ \Omega(\sum_{i=1}^r \xi_i \vP_i) = \Pi \begin{pmatrix} * & \Xi \\
  -\Xi & 0 \end{pmatrix} \Pi^T \] where
\begin{equation}
\label{xi}
\Xi = \begin{pmatrix} 
  (n-2) \xi_1^2 & -2 \xi_1 \xi_2 & \ldots & -2 \xi_1 \xi_{r} \\
  -2 \xi_1 \xi_2 & (n-2) \xi_2^2 & \ldots & -2 \xi_2 \xi_r \\
  \vdots & \vdots & \ddots & \vdots \\
  -2 \xi_1 \xi_{r} & -2 \xi_2 \xi_r & \ldots & (n-2) \xi_{r}^2 \end{pmatrix} 
\end{equation}
and $\Pi$ is the $n \times 2r$ matrix with columns $\vP_1, \ldots,
\vP_r, d\vP_1, \ldots, d\vP_r$.
\end{Lemma}
\begin{Proof}
  Recall that $\Omega = n B - 2 A$. The case $r=2$ is immediate from
  Lemma~\ref{lem:5}. Since the entries of $\Omega$ are quadratic forms
  the general case follows.
\end{Proof}

We now check the hypotheses {\bf (H1)}, {\bf (H2)} and {\bf (H3)}.

\medskip

\begin{ProofOf}{(H1) and (H3)}
Suppose $n-2r \ge 1$. A generic point $P \in \Sec^r C$ may be written
$P = [\sum_{i=1}^r \xi_i \vP_i ]$ for some $0_E \not= P_1, \ldots, P_r
\in E$ distinct and $\xi_1, \ldots, \xi_r \not= 0$. 
By Proposition~\ref{lem:tangent} the tangent space $T_P \Sec^r C
\subset \PP^{n-1}$ is spanned by $\vP_1, \ldots, \vP_r, d \vP_1,
\ldots d \vP_r$. In particular these $2r$ vectors are linearly
independent.

For $f \in I(\Sec^r C)$ we have $\sum_{i=1}^n \frac{\partial
  f}{\partial x_i}(P) \vv_i = 0$ for any $\vv$ in the linear span of
$\vP_1, \ldots, \vP_r, d \vP_1, \ldots d \vP_r$. By Lemma~\ref{lem:6}
the columns of $\Omega$ are linear combinations of these vectors. So
for each $1 \le j \le n$ the form $\sum_{i=1}^n \frac{\partial
  f}{\partial x_i} \Omega_{ij}$ vanishes at $P$.  Since $P \in \Sec^r
C$ is generic, 
this proves {\bf (H1)}. Since $n \notin \{0, 2r\}$ and 
$\xi_1, \ldots, \xi_r \not= 0$, the matrix~\eqref{xi} 
is non-singular. 
Therefore $\rank \Omega(P) = 2r$ and this proves {\bf (H3)}. 
\end{ProofOf}

\begin{ProofOf}{(H2)}
We write $n=2r$ and $\Sec^{r-1} C = \{F_1 = F_2 =
0\}$ where $F_1$ and $F_2$ are forms of degree $r$. We must show that
the form
\begin{equation}
\label{mixed}
\sum_{i,j=1}^n \frac{\partial F_1}{\partial x_i} \Omega_{ij} \frac{\partial F_2}{\partial x_j}
\end{equation}
is identically zero. A generic point $P \in \Sec^r C = \PP^{n-1}$ may
be written $P = [\sum_{i=1}^r \xi_i \vP_i ]$ for some $0_E \not= P_1,
\ldots, P_r \in E$ distinct and $\xi_1, \ldots, \xi_r \not= 0$. In
addition we may assume that $2(P_1 + \ldots + P_r) \not\sim H$ where
$H$ is the hyperplane section. This ensures that the vectors $\vP_1,
\ldots, \vP_r, d \vP_1, \ldots d \vP_r$ are linearly independent. We
choose co-ordinates on $\PP^{n-1}$ so that $[\vP_1] = (1:0:
\ldots:0)$, $[\vP_2] = (0:1:0: \ldots:0)$, \ldots,
$d \vP_r =(0: \ldots :0:1)$. Since $F_1$ and $F_2$
vanish on $\Sec^{r-1}C$ they vanish on the linear span of any $r-1$
of the $[\vP_i]$.  Replacing $F_1$ and $F_2$ by suitable linear
combinations we may assume
\begin{align*}
  F_1(x_1, \ldots, x_r, 0, \ldots, 0) & = 0, \\
  F_2(x_1, \ldots, x_r, 0, \ldots, 0) & = x_1 x_2 \ldots x_r.
\end{align*}
Therefore at $P = (\xi_1: \ldots : \xi_r:0: \ldots: 0)$ we have
\begin{align*}
  (\tfrac{\partial F_1}{\partial x_1}(P), \ldots, \tfrac{\partial
    F_1}{\partial x_n}(P))
  & = (0, \ldots, 0, *, \ldots, *), \\
  (\tfrac{\partial F_2}{\partial x_1}(P), \ldots, \tfrac{\partial
    F_2}{\partial x_n}(P)) & = (\prod_{i \not=1} \xi_i, \ldots,
  \prod_{i\not= r} \xi_i, *, \ldots, *).
\end{align*}
By Lemma~\ref{lem:6} we have
\[ \Omega(P) = \begin{pmatrix} * & \Xi \\ -\Xi & 0 \end{pmatrix} \]
where $\Xi$ is the matrix~\eqref{xi}.  Since $n=2r$ the coefficients
in each row and column of $\Xi$ sum to zero.  Therefore the
form~\eqref{mixed} vanishes at $P$.  Since $P \in \PP^{n-1}$ is
generic, this shows that the form is identically zero.
\end{ProofOf}


This completes the proof of Theorem~\ref{thm1}.

\section{Explicit formulae}
\label{sec:pf2}

In this section we give an explicit formula for the matrix $\Omega$
defined in Section~\ref{sec:pf1}.  As before $E$ is the elliptic curve
\begin{equation*}
\label{weierE2}
y^2 + a_1 xy + a_3 y = x^3 + a_2 x^2 + a_4 x + a_6.
\end{equation*}
with invariant differential $\omega = dx/(2y + a_1 x + a_3)$.  We
embed $E$ in $\PP^{n-1}$ via
\begin{equation}
\label{emb}
(x_0: x_2: x_3 : \ldots : x_n )  = (1,x,y,x^2,xy,x^3,x^2y,x^4, \ldots ) 
\end{equation}
Notice there is no $x_1$. The indicator function of a set $X$ is
denoted $\one_X$. We define linear forms in indeterminates $\{x_m : m
\in \Z\}$ as follows
\begin{align*}
  \dot{x}_m &= \tfrac{m}{2}(2 x_{m+1} + a_1 x_m + a_3 x_{m-2}) +
  \one_{\text{$m$ odd}} \sum_{i=1}^6
  (-1)^i (m - \tfrac{i}{2}) a_i x_{m+1-i} \\
  \overline{x}_m &= \tfrac{1}{2}(2 x_{m+1} + a_1 x_{m} + a_3 x_{m-2})
  + \one_{\text{$m$ odd}} \sum_{i=1}^6 (-1)^i a_i x_{m+1-i}
\end{align*}
where by convention $a_5= 0$. For $x \in \R$ we let $\sign(x)=-1,0,1$
according as $x$ is negative, zero or positive. For $r,s \in \Z$ we
define
\begin{align*}
  A_{rs} &= x_r \dot{x}_s - x_s \dot{x}_r,  \\
  B_{rs} &= \sum_{k =-\infty}^{\infty} \sign(k+\tfrac{1}{2}) (x_{r+2k}
  \overline{x}_{s-2k} - x_{s+2k} \overline{x}_{r-2k}).
\end{align*}

\begin{Theorem}
  \label{thm:explicitomega}
Let $C \subset \PP^{n-1}$ be the image of $E$ under the embedding~\eqref{emb}.
\begin{enumerate}
\item $A = (A_{rs})_{r,s = 0,2,\ldots,n}$ and $B= (B_{rs})_{r,s =
    0,2,\ldots,n}$ are $n \times n$ alternating matrices of quadratic
  forms in $x_0,x_2, \ldots, x_{n+1}$.
\item $\Omega = n B - 2 A$ is an $n \times n$ alternating matrix of
  quadratic forms in $x_0,x_2, \ldots,$ $x_n$. It satisfies the
  conclusions of Theorem~\ref{thm1} and
  \begin{equation}
\label{eq:explicit}
(n-2) \omega = \frac{ x_j^2 d(x_i/x_j) }{\Omega_{ij}(x_1, \ldots,x_n)} \qquad
\text{ for all } i \not= j. 
\end{equation}
\end{enumerate}
\end{Theorem}

\begin{Proof}
  It is part of the theorem that the indeterminates $x_m$ for $m
  \notin \{0,2,3, \ldots, n\}$ cancel from the formula for
  $\Omega$.  So when applying the theorem we simply set them to be
  zero. However we will not do this in the proof. Since
  $\overline{x}_m$ is a linear combination of $x_{m+1}, x_m, \ldots,
  x_{m-5}$ each $B_{rs}$ is of the form $\sum_{ij} c_{ij} x_i x_j$
  where each $c_{ij}$ is a finite sum. But it is not immediately clear
  that the $B_{rs}$ are polynomials, i.e. that $c_{ij} = 0$ for all
  but finitely many pairs $(i,j)$.  We check this first.

  If $r \equiv s \pmod{2}$ and $r<s$ then
  \begin{equation}
  \label{beven}
  B_{rs} = 2(x_r \overline{x}_s + x_{r+2} \overline{x}_{s-2} + \ldots + x_{s-2} 
  \overline{x}_{r+2}) 
\end{equation}
whereas if $r$ is even and $s$ is odd then
\begin{equation}
\label{bodd}
B_{rs} = - a_1 x_r x_s +
Q_{r,s+1} + a_2 Q_{r,s-1} + a_4 Q_{r,s-3} + a_6 Q_{r,s-5} - Q_{s,r+1} 
\end{equation}
where 
\[ Q_{ij} = \left\{ \begin{array}{ll}
    x_i x_j + x_{i+2} x_{j-2} + \ldots + x_j x_i & \text{ if } i < j+2, \\
    0 & \text{ if } i=j+2, \\
    -(x_{i-2} x_{j+2} + x_{i-4} x_{j+4} + \ldots + x_{j+2} x_{i-2} ) &
    \text{ if } i > j+2. 
  \end{array} \right. \]
Since $B_{sr} = -B_{rs}$ this proves that the $B_{rs}$
are polynomials.

We show that the matrices $A$ and $B$ defined in the statement of the
theorem represent the linear maps $A$ and $B$ defined in
Section~\ref{sec:pf1}. The theorem then follows from the results of
Sections~\ref{sec:outline}, \ref{sec:pf1} and~\ref{sec:checkhyp}.  In
particular~\eqref{eq:explicit} follows from Lemma~\ref{lem:omega}.

In the statement of the theorem the $\{x_m : m \in \Z\}$ are
indeterminates.  However for the proof they will be the following
rational functions on $E$,
\[ x_m = \left\{ \begin{array}{ll} x^{m/2} & \text{ if $m$ is even, } \\
    x^{(m-3)/2} y & \text{ if $m$ is odd.} \end{array} \right. \] As
rational functions on $E$, we claim that $\dot{x}_m = dx_m/\omega$ (in
agreement with the notation in Section~\ref{sec:pf1}) and
$\overline{x}_m = \tfrac{1}{2} x_{m-2} (2 y + a_1 x + a_3)$. In
checking these claims, we start with the right hand sides, since this
also serves to motivate the definitions of $\dot{x}_m$ and
$\overline{x}_m$.  For $m$ even we have
\begin{align*}
  dx_m/\omega &= \tfrac{m}{2} x^{(m-2)/2} (dx/\omega) \\
  & = \tfrac{m}{2} x^{(m-2)/2} (2 y + a_1 x + a_3) \\
  & =\tfrac{m}{2} (2 x_{m+1} + a_1 x_m + a_3 x_{m-2}),
\end{align*}
and
\[ \tfrac{1}{2} x_{m-2} (2 y + a_1 x + a_3) =\tfrac{1}{2} (2 x_{m+1} +
a_1 x_m + a_3 x_{m-2}). \] For $m$ odd we have
\begin{align*}
  dx_m/\omega &= \tfrac{m-3}{2}x^{(m-5)/2}y(dx/\omega)+x^{(m-3)/2}(dy/\omega) \\
  &= \tfrac{m-3}{2}x^{(m-5)/2}(2 y^2 + a_1 xy + a_3y)
  + x^{(m-3)/2} (3x^2 + 2a_2 x + a_4 - a_1 y) \\
  &= \tfrac{m-3}{2}x^{(m-5)/2}(-a_1 xy-a_3 y+2x^3 + 2a_2 x^2 + 2a_4 x
  + 2a_6)
  \\ & \hspace{12em} + x^{(m-5)/2} (3x^3 + 2a_2 x^2 + a_4x - a_1 xy) \\
  &= x^{(m-5)/2}\big(m x^3 - \tfrac{m-1}{2} a_1 xy- \tfrac{m-3}{2} a_3
  y + \textstyle\sum_{i=1}^3 (m-i) a_{2i} x^{3-i}\big) \\
  &= \tfrac{m}{2}(2 x_{m+1} + a_1 x_{m} + a_3 x_{m-2}) +
  \textstyle\sum_{i=1}^6 (-1)^i (m - \tfrac{i}{2}) a_i x_{m+1-i},
\end{align*}
and
\begin{align*}
  \tfrac{1}{2} x_{m-2} (2 y + a_1 x + a_3)
  &=  \tfrac{1}{2}x^{(m-5)/2}(2 y^2 + a_1 xy + a_3y) \\
  &=  \tfrac{1}{2}x^{(m-5)/2} (-a_1 xy-a_3 y+2x^3 + 2a_2 x^2 + 2a_4 x + 2a_6) \\
  &= \tfrac{1}{2}(2 x_{m+1} + a_1 x_{m} + a_3 x_{m-2}) +
  \textstyle\sum_{i=1}^6 (-1)^i a_i x_{m+1-i}.
\end{align*}

It is now clear that $A(x_r \wedge x_s) = A_{rs}$ for all $r,s \in
\Z$.  It remains to prove the same for $B$. By definition of $B$ we
have
\[ B(x_r \wedge x_s) = \frac{y_P + y_Q + a_1 x_Q + a_3}{x_P - x_Q}
(x_r(Q) x_s(P) - x_r(P) x_s(Q)) \bigg|_{P=Q}
\]
where $P, Q$ are points on $E$.  Since $\overline{x}_m = \tfrac{1}{2}
x_{m-2} (2 y + a_1 x + a_3)$ we have
\begin{align*}
  2 x_r(P) \overline{x}_s(Q)
  &= (2 y_Q + a_1 x_Q + a_3) x_r(P) x_{s-2}(Q) \\
  &= \frac{2 y_Q + a_1 x_Q + a_3}{x_P - x_Q} (x_{r+2}(P) x_{s-2}(Q) -
  x_r(P) x_s(Q)).
\end{align*}
Adding this to the same expression with $(r,s)$ replaced by
$(s-2,r+2)$ and then setting $P=Q$ gives
\begin{equation} 
\label{myid}
B_{rs} - B_{r+2,s-2} = 2 (x_r \overline{x}_s + x_{s-2} \overline{x}_{r+2}) = 
B(x_r \wedge x_s) - B(x_{r+2} \wedge x_{s-2}). 
\end{equation}
Rather more obviously, replacing $(r,s)$ by $(r+2,s+2)$ changes
$B_{rs}$ and $B(x_r \wedge x_s)$ in the same way, that is, by shifting
the subscripts up by $2$. So to prove $B(x_r \wedge x_s) = B_{rs}$ for
all $r,s \in \Z$ it suffices to prove it for all $r \in \{0,1\}$ and
$s \in \{0,1,2,3\}$. This is a finite calculation. We give two
examples:
\begin{align*}
  B(x_0 \wedge x_3) &= \frac{y_P + y_Q + a_1 x_Q + a_3}{x_P - x_Q}
  (y_P - y_Q) \bigg|_{P=Q} \\
  &= \frac{(y_P^2 + a_1 x_P y_P + a_3 y_P)
    - (y_Q^2 + a_1 x_Q y_Q + a_3 y_Q)}{x_P - x_Q} - a_1 y_P \bigg|_{P=Q} \\
  &= (x_P^2 + x_P x_Q + x_Q^2 - a_1 y_P + a_2(x_P + x_Q) + a_4) \big|_{P=Q} \\
  &= 2x_0 x_4 + x_2^2 - a_1 x_0 x_3 + 2 a_2 x_0 x_2 + a_4 x_0^2,
\end{align*}
and
\begin{align*}
  B(x_2 \wedge x_3) &= \frac{y_P + y_Q + a_1 x_Q + a_3}{x_P - x_Q}
  (y_P(x_Q-x_P) + x_P(y_P - y_Q)) \bigg|_{P=Q} \\
  & = (-y_P (y_P + y_Q + a_1 x_Q + a_3) + x_P(x_P^2 + x_P x_Q +
  \ldots 
  + a_4)) \big|_{P=Q} \\
  &= (x_P^2x_Q + x_P x_Q^2 - y_P y_Q - a_1 x_Q y_P
  + a_2 x_P x_Q - a_6) \big|_{P=Q} \\
  &= 2x_2 x_4 - x_3^2 - a_1 x_2 x_3 + a_2 x_2^2 - a_6 x_0^2.
\end{align*}
It is easy to check using~\eqref{bodd} that these are equal to
$B_{03}$ and $B_{23}$. The other cases we need can then be checked
using~\eqref{myid} and the fact that $B$ is alternating.
\end{Proof}

\section{Proof of Theorem~\ref{thm2}}
\label{sec:scale}

Let $\Omega = nB - 2A$ be as in Theorem~\ref{thm:explicitomega}.  Then
$c_4(\Omega) = f_n(a_1, \ldots,a_6)$ and $c_6(\Omega) = g_n(a_1,
\ldots,a_6)$ for some polynomials $f_n$ and $g_n$.  We consider the
effect of a change of Weierstrass equation, with notation as in
\cite[Chapter III]{Sil}.

\begin{Lemma}
\label{lem:urst}
Let $a_1, \ldots, a_6$ and $a'_1, \ldots, a'_6$ be the coefficients of
two Weierstrass equations related by $x = u^2 x' + r$ and $y = u^3 y'
+ u^2 s x' + t$. Then
\begin{align*}
  f_n(a_1, \ldots, a_6) &= u^{4} f_n(a'_1, \ldots, a'_6) \\
  g_n(a_1, \ldots, a_6) &= u^{6} g_n(a'_1, \ldots, a'_6)
\end{align*}
\end{Lemma}
\begin{Proof}
  This follows from Corollary~\ref{howc4c6change} and $u^{-1} \omega'
  = \omega$.
\end{Proof}

It follows by Lemma~\ref{lem:urst}, and the standard procedure for
converting a Weierstrass equation to the shorter form $y^2 = x^3 + a x
+ b$, that $f_n$ and $g_n$ are scalar multiples of the usual
polynomials $c_4$ and $c_6$ in $a_1, \ldots, a_6$.  Explicitly,
\begin{equation}
\label{w-inv}
\begin{aligned}
  f_n(a_1, \ldots, a_6) &= \xi_n (b_2^2 - 24 b_4) = \xi_n(a_1^4 + \ldots), \\
  g_n(a_1, \ldots, a_6) &= \eta_n (-b_2^3 + 36 b_2 b_4 - 216 b_6) =
  \eta_n( -a_1^6 + \ldots),
\end{aligned}
\end{equation}
where $b_2 = a_1^2 + 4a_2$, $b_4 = 2 a_4 + a_1 a_3$ and $b_6 = a_3^2 +
4a_6$.

To complete the proof of Theorem~\ref{thm2} we must compute the
constants $\xi_n$ and $\eta_n$.  For any given value of $n$ these can
be read off from a single numerical example. However we need to
compute these constants for all $n$.  We write
\[\Omega = \Omega^{(0)} + a_1 \Omega^{(1)} + a_2 \Omega^{(2)} 
+ a_3 \Omega^{(3)} + a_4 \Omega^{(4)} + a_6 \Omega^{(6)}. \] Since
$c_4(\Omega)$ and $c_6(\Omega)$ have degrees $4$ and $6$ in the
coefficients of the entries of $\Omega$, we see by~\eqref{w-inv} that
it suffices to compute the invariants of $\Omega^{(1)}$.

We put
\[ \gamma_{rs} = (-1)^{\max(r,s)} \sign(s-r) n \,\, - \,\, 2\left( (-1)^s \lfloor 
\tfrac{s}{2} \rfloor - (-1)^r \lfloor \tfrac{r}{2} \rfloor \right). \]
\begin{Lemma}
\label{getom1}
The alternating matrix $\Omega^{(1)}$ has entries above the diagonal
\begin{equation}
\label{om1}
 \gamma_{rs} x_r x_s + (-1)^s n \one_{r+s \text{ \rm{even}}} 
\textstyle\sum_{k=1}^{(s-r)/2-1} x_{r+2k} x_{s-2k}.
\end{equation}
\end{Lemma}
\begin{Proof} Since $\Omega = n B - 2 A$ we have $\Omega^{(1)} = n
  B^{(1)} - 2 A^{(1)}$ where the superscripts indicate that we are
  taking the coefficient of $a_1$. Then $A^{(1)}$ has $(r,s)$ entry
  \[ \left( (-1)^s \lfloor \tfrac{s}{2} \rfloor - (-1)^r \lfloor
    \tfrac{r}{2} \rfloor \right) x_r x_s \] whereas~\eqref{beven}
  and~\eqref{bodd} show that if $r < s$ then $B^{(1)}$ has $(r,s)$
  entry
  \[ \left\{ \begin{array}{ll} (-1)^s (x_r x_s + x_{r+2} x_{s-2} +
      \ldots x_{s-2} x_{r+2})
      & \text{ if } r \equiv s \pmod{2} \\
      (-1)^s x_r x_s & \text{ if } r \not\equiv s \pmod{2}.
    \end{array} \right. \]
\end{Proof}

\begin{Lemma}
  The matrices $\Omega^{(1)}$, $\Omega' = \big(\gamma_{rs} x_r
  x_s\big)_{r,s = 0,2,3, \ldots,n}$ and
  \[ \Lambda = \big((\sign(j-i) n - 2(j-i)) x_i x_j\big)_{i,j =0,1,
    \ldots,n-1}\] all have the same invariants $c_4$ and $c_6$.
\end{Lemma}

\begin{Proof}
  We first explain why $\Omega^{(1)}$ and $\Omega'$ have the same
  invariants, despite the ``extra terms'' in~\eqref{om1}.  We start
  with $\Omega^{(1)}$.  The only entries involving $x_0$ are in the
  first row and column.  We replace $x_0$ by $\lambda^{-1} x_0$ and
  multiply the first row and column by $\lambda$. By
  Lemma~\ref{leminv} this does not change the invariants, but setting
  $\lambda = 0$ removes the extra terms from the first row and
  column. Now the only entries involving $x_2$ are in the second row
  and column. We replace $x_2$ by $\lambda^{-1} x_2$ and multiply the
  second row and column by $\lambda$. This does not change the
  invariants, but setting $\lambda = 0$ removes the extra terms from
  the second row and column. We repeat this procedure for all
  subsequent rows and columns.  In the end we remove all the extra
  terms, and are left with the matrix $\Omega'$.

  We define a bijection $\pi: \{0,1, \ldots, n-1\} \to
  \{0,2,3,\ldots,n\}$ by
  \[ \pi(i) = \left\{ \begin{array}{ll} 2i & \text{ if } i \le n/2, \\
      2(n-i) + 1 & \text{ if } i > n/2. \end{array} \right. \] We then
  compute
  \[ \gamma_{\pi(i),\pi(j)} = \left\{ \begin{array}{ll} \sign(j-i) n -
      2(j-i)
      & \text{ if $i \le n/2$ and $j \le n/2$,} \\
      -n - 2(-(n-j) - i)
      & \text{ if $i \le n/2$ and $j > n/2$, } \\
      n - 2( j + (n-i))
      & \text{ if $i > n/2$ and $j \le n/2$,} \\
      \sign(j-i) n - 2(-(n-j)+(n-i)) & \text{ if $i > n/2$ and $j >
        n/2$.} \end{array} \right. \] In all cases we have
  $\gamma_{\pi(i),\pi(j)} = \sign(j-i) n - 2(j-i)$.  Therefore
  $\Omega'$ and $\Lambda$ are related by a permutation matrix. It
  follows by Lemma~\ref{leminv} that they have the same invariants.
\end{Proof}

\begin{Lemma}
The alternating matrix of quadratic forms
\[ \Lambda = \begin{pmatrix}
  0 & (n-2) x_1 x_2 & (n-4) x_1 x_3 & (n-6) x_1 x_4 & \cdots & (2-n) x_1 x_n \\
  & 0 & (n-2) x_2 x_3 & (n-4) x_2 x_4 & \cdots & (4-n) x_2 x_n \\
  &  & 0 & (n-2) x_3 x_4 & \cdots & (6-n) x_3 x_n \\
  & - & & && \vdots \\
  & &  &   \ddots &  & (n-2) x_{n-1} x_{n} \\
  & & & & & 0
\end{pmatrix}
\] 
has invariants $c_4(\Lambda) = (n-2)^4$ and $c_6(\Lambda) = -(n-2)^6$.
\end{Lemma}

\begin{Proof} We have $\Lambda = (\lambda_{rs} x_r x_s)_{r,s = 1,
    \ldots,n}$ where $\lambda_{rs} = \sign(s-r) n - 2(s-r)$.
  Following the definitions of $c_4$ and $c_6$ in
  Section~\ref{sec:stat} we put
  \begin{align*}
    M_{ij} &= \sum_{r,s=1}^n \frac{\partial \Lambda_{ir}}{\partial
      x_s} \frac{ \partial \Lambda_{js}}{\partial x_r}
    = \mu_{ij} x_i x_j \\
    N_{ijk} &= \sum_{r=1}^n \frac{\partial M_{ij}}{\partial x_r}
    \Lambda_{rk} = \nu_{ijk} x_i x_j x_k
\end{align*}
where $\mu_{ij} = (\textstyle\sum_{r=1}^n \lambda_{ir} \lambda_{jr}) -
\lambda_{ij}^2$ and $\nu_{ijk} = \mu_{ij} (\lambda_{ik} +
\lambda_{jk})$.  It is not hard to show that
\begin{align*}
  \sum_{r=1}^n \sign(i-r) \sign(j-r) &= n - 2|i-j| - \delta_{ij}, \\
  \sum_{r=1}^n (i-r) \sign(j-r) &= 2ij - j^2 - (n+1)i + n(n+1)/2, \\
  \sum_{r=1}^n (i-r) (j-r) &= nij - (i+j) n(n+1)/2 + n(n+1)(2n+1)/6.
\end{align*}
We use these to compute 
\[ \sum_{r=1}^n \lambda_{ir} \lambda_{jr} = 2n|i-j|^2 - 2n^2|i-j| -
\delta_{ij} n^2 + (n^3 + 2n)/3 \] and then subtract off
\[\lambda_{ij}^2= 4 |i-j|^2 - 4n |i-j| + (1- \delta_{ij}) n^2\]
to get
\[ \mu_{ij} = 2 (n-2) ( |i-j|^2 - n|i-j|) + n(n-1)(n-2)/3. \] Noting
the symmetries $\mu_{ij} = \mu_{ji}$ and $\nu_{ijk} = \nu_{jik}$, and
using computer algebra to check our calculations, we find
\[ \sum_{i,j,r,s=1}^n \frac{\partial^2 M_{ij}}{\partial x_r \partial
  x_s} \frac{\partial^2 M_{rs}}{\partial x_i \partial x_j} = 4 \sum_{i
  \le j} \mu_{ij}^2 = (16/3) n (n-2)^2 \binom{n+3}{5} \] and
\begin{align*}
  \sum_{i,j,k,r,s,t=1}^n & \frac{\partial^3 N_{ijk}}{\partial
    x_r \partial x_s \partial x_t}
  \frac{\partial^3 N_{rst}}{\partial x_i \partial x_j \partial x_k} \\
  &= 4 \sum_{i \le j \le k} (\nu_{ijk}+\nu_{jki}+\nu_{kij})^2 \\
  &= 4 \sum_{i \le j \le k} (\la_{ij}(\mu_{ik} - \mu_{jk}) 
  + \la_{jk}(\mu_{ij} - \mu_{ik}) + \la_{ik} (\mu_{ij} - \mu_{jk}))^2 \\
  & = 64 (n-2)^2 \sum_{i \le j \le k} (i-2j+k)^2 (n+i+j-2k)^2 (n+2i-j-k)^2 \\
  & = 64 n (n-2)^3 \binom{n+5}{7}.
\end{align*}
The final sums are evaluated using the standard formulae for
$\sum_{i=1}^n i$, $\sum_{i=1}^n i^2$, etc. In practice it is simpler
to observe that the answer is a polynomial in $n$, say of degree at
most $d$, and then check the result for $d+1$ distinct values of $n$.

Finally scaling by the constants included in the
definitions~\eqref{c4} and \eqref{c6} it follows that $c_4(\Lambda) =
(n-2)^4$ and $c_6(\Lambda) = -(n-2)^6$.
\end{Proof}

The last two lemmas show that $\xi_n = (n-2)^4$ and $\eta_n =
(n-2)^6$. Therefore $c_4(\Omega) = (n-2)^4 c_4(E)$ and $c_6(\Omega) =
(n-2)^6 c_6(E)$. Let $\omega = dx/(2y + a_1x + a_3)$. By the 
formulae in \cite[Chapter III]{Sil} we have
\[ (E, \omega) \cong (y^2 = x^3 - 27 c_4(E) x - 54
c_6(E), 3 dx/y ). \] 
Therefore 
\[ (E,(n-2) \omega) \cong (y^2 = x^3 - 27 c_4(\Omega) x - 54
c_6(\Omega), 3 dx/y ). \] 
Recalling from Theorem~\ref{thm:explicitomega} 
that $\Omega = n B - 2 A$ represents the invariant 
differential $(n-2) \omega$, this completes the proof of
Theorem~\ref{thm2}.

\section{Higher secant varieties}
\label{sec:hsec}

In this final section we give references and proofs for the facts
about higher secant varieties we used earlier in the paper.

\begin{Theorem} 
\label{thm:hsec-eqns}
Let $C \subset \PP^{n-1}$ be a genus one normal curve of degree $n \ge
3$.
\begin{enumerate}
\item $\Sec^r C \subset \PP^{n-1}$ is an irreducible variety of 
codimension $\max(n-2r,0)$. 
\item The vector space of forms of degree $r+1$ vanishing on $\Sec^r
  C$ has dimension $\beta(r+1,n)$, where
\begin{equation*}
  \beta(r,n) = \binom{n-r}{r} + \binom{n-r-1}{r-1}
\end{equation*}
is the number of ways of choosing $r$ elements from $\Z/n\Z$ such that
no two elements are adjacent.
\item If $n-2r \ge 2$ then the homogeneous ideal $I(\Sec^r C)$ is
  generated by forms of degree $r+1$.
\item if $n-2r = 1$ then $\Sec^r C$ is a hypersurface of degree $n$.
\item If $n-2r \ge 1$ then $\Sec^r C$ has singular locus
  $\Sec^{r-1}C$.
\end{enumerate}
\end{Theorem}
\begin{Proof} (i) This is a general fact about curves.
  See for example \cite[Section 1]{Lange}. \\
  (ii), (iii), (iv). More generally the minimal free resolution for
  $I(\Sec^r C)$ was computed in \cite[Section 8]{vBH}.  See
  \cite[Section 5]{GP} for the cases $r=1,2$,
  and~\cite[Section 4]{pfpres} for further discussion. \\
  (v) This is \cite[Proposition 8.15]{vBH}.
\end{Proof}

\subsection{Computing equations for higher secant varieties}

The following two propositions may be used to compute equations for
$\Sec^r C$ from equations for $C$.  We say that a form $f$ vanishes on
$C$ with multiplicity $r$ if (passing to affine co-ordinates) the
Taylor expansion of $f$ at each point $P \in C$ begins with terms of
order greater than or equal to $r$.

\begin{Proposition} 
\label{lem:deriv}
Let $C \subset \PP^{n-1}$ be a variety contained in no hyperplane.
Let $f$ be a form of degree $r+1$.
\begin{enumerate}
\item If $r \ge 1$ then
\[ f \in I(\Sec^r C) \iff \text{$f$ vanishes on $C$ with multiplicity $r$.} \]
\item If $r \ge 2$ then 
  \[ f \in I(\Sec^r C) \iff \frac{\partial f}{\partial x_i} \in
  I(\Sec^{r-1} C) \text{ for all } i=1,\ldots, n.\]
\end{enumerate}
\end{Proposition}
\begin{Proof}
  (i) We choose $P_1, \ldots, P_n \in C$ spanning $\PP^{n-1}$. By a
  change of co-ordinates we may assume $P_1=(1:0: \ldots :0), P_2 =
  (0:1:0 : \ldots:0), \ldots , P_n=(0:0: \ldots :1)$. If $f \in
  I(\Sec^r C)$ then it vanishes on the linear span of any $r$ of the
  $P_i$.  Therefore the monomials appearing in $f$ involve at least
  $r+1$ of the $x_i$, and since $f$ has degree $r+1$ must be
  squarefree.  But then $f$ vanishes at $P_1$ with multiplicity
  $r$. Since $P_1 \in C$ was arbitrary it follows that $f$ vanishes on
  $C$ with multiplicity $r$.
 
  Conversely, suppose $f$ vanishes on $C$ with multiplicity $r$.  Let
  $\Pi$ be an $(r-1)$-plane spanned by points $P_1, \ldots,P_r \in C$.
  By a change of co-ordinates we may assume $P_1=(1:0: \ldots :0), P_2
  = (0:1:0 : \ldots:0), \ldots $.  Then $f(x_1, \ldots, x_r,0, \ldots,
  0)$ has total degree $r+1$, but has degree at most $1$ in each of
  the variables.  It follows that $f$ vanishes on $\Pi$. By definition
  $\Sec^r C$ is the Zariski closure of the union of all such
  $(r-1)$-planes. Therefore $f \in I(\Sec^r C)$ as required.  \\
  (ii) Since $\Char(K) = 0$ this follows from (i).
\end{Proof}

Now let $C \subset \PP^{n-1}$ be a genus one normal curve.
Taking $r=1$ in Theorem~\ref{thm:hsec-eqns} shows that the homogeneous
ideal $I(C)$ is generated by a vector space of quadrics of dimension
$n(n-3)/2$.  Suppose we know a basis for this space. Then by
repeatedly applying Proposition~\ref{lem:deriv}(ii) we can find a
basis for the space of forms of degree $r+1$ vanishing on $\Sec^r C$.
Theorem~\ref{thm:hsec-eqns}(iii) tells us that if $n - 2r \ge 2$ then
these forms define $\Sec^r C$. The following proposition covers the
remaining case.

\begin{Proposition}
\label{lem:deriv1}
Suppose $n-2r = 1$. Let $f$ be a form of degree $n$. If $r \ge 2$ then
\[ f \in I(\Sec^r C) \iff \frac{\partial f}{\partial x_i} \in
I(\Sec^{r-1} C)^2 \text{ for all } i=1,\ldots, n.\]
\end{Proposition}
\begin{Proof}
  ``$\Rightarrow$'' Let $H$ be the divisor of a hyperplane section,
  and let $P \in C$ be any point. Let $C_+ \subset \PP^n$ and $C_-
  \subset \PP^{n-2}$ be the images of $C$ embedded via the linear
  systems $|H \pm P|$.  We choose co-ordinates so that the
  isomorphisms $C_+ \to C \to C_-$ are given by
  \[ (x_1: \ldots : x_{n+1}) \mapsto (x_1: \ldots : x_n) \mapsto (x_1
  : \ldots : x_{n-1}). \] In particular $P$ is the point $(x_1: \ldots
  : x_n) = (0: \ldots :0 :1)$.  By Theorem~\ref{thm:hsec-eqns} we know
  that $I(\Sec^{r-1} C_-)$ is generated by forms $g_1,g_2 \in K[x_1,
  \ldots,x_{n-1}]$ of degree $r$.  By \cite[Corollary 2.3]{pfpres}
  there exist forms $h_1,h_2 \in K[x_1,\ldots,x_n]$ of degree $r+1$
  such that $f_i = x_{n+1} g_i + h_i \in I(\Sec^r C_+)$ for $i=1,2$.
  Then $F = g_1 h_2 - g_2 h_1$ belongs to
  \[ I(\Sec^r C_+) \cap K[x_1,\ldots,x_n] = I(\Sec^r C). \] Since
  $g_1,g_2$ are coprime and $f_1,f_2$ are irreducible it is clear that
  $F$ is non-zero. By Theorem~\ref{thm:hsec-eqns}(iv) we have
  $I(\Sec^rC ) = (F)$. We compute
  \[ \frac{\partial F}{\partial x_n} = \frac{\partial f_1}{\partial
    x_{n+1}} \, \frac{\partial f_2}{\partial x_n} - \frac{\partial
    f_1}{\partial x_n} \, \frac{\partial f_2}{\partial x_{n+1}}.  \]
  On the other hand, for $i=1,2$ and $j = n,n+1$ we have
  \[ \frac{\partial f_i}{\partial x_j} \in I(\Sec^{r-1} C_+) \cap
  K[x_1,\ldots,x_n] = I(\Sec^{r-1} C). \] Therefore $\frac{\partial
    F}{\partial x_n} \in I(\Sec^{r-1} C)^2$.
  Since $P \in C$ was arbitrary, and $C$ spans $\PP^{n-1}$, 
  the result follows. \\
  ``$\Leftarrow$'' Let $P_1, \ldots, P_r$ be $r$ distinct points on
  $C$.  By a change of co-ordinates we may assume $P_1=(1:0: \ldots
  :0), P_2 = (0:1:0 : \ldots:0), \ldots $.  By
  Proposition~\ref{lem:deriv} we know that $f$ vanishes on $C$ with
  multiplicity $2(r-1) + 1 = n-2$.  Therefore $f(x_1, \ldots, x_r,0,
  \ldots, 0)$ has total degree $n$, but has degree at most $2$ in each
  of the variables. Since $2 r < n$ it follows that $f$ vanishes on
  the linear span of $P_1, \ldots, P_r$.  By definition $\Sec^r C$ is
  the Zariski closure of the union of all such
  $(r-1)$-planes. Therefore $f \in I(\Sec^r C)$ as required.
\end{Proof} 

\subsection{Proof of Proposition~\ref{lem:tangent}}
\label{sec:tgt}

Let $C \subset \PP^{n-1}$ be a genus one normal curve of degree
$n$. Let $H$ be the divisor of a hyperplane section. We identify
$\LL(H)$ with the space of linear forms on $\PP^{n-1}$. For $D$ an
effective divisor on $C$ we write $\overline{D} \subset \PP^{n-1}$ for
the linear subspace cut out by $\LL(H-D) \subset \LL(H)$. We have
\[\Sec^r C = \bigcup_{\deg D = r} \Dbar.\] 
We also put $D^\circ = \Dbar \setminus \cup_{D'<D} 
\overline{D'}$. The $\gcd$ and $\lcm$
of divisors $\sum m_P P$ and $\sum m'_P P$ are $\sum \min(m_P,m'_P) P$
and $\sum \max(m_P,m'_P) P$.

\begin{Lemma}
\label{Dbarlemma}
Let $D$, $D_1$, $D_2$ be effective divisors on $C$.  \\
(i) If $\deg D < n$ then $\dim \Dbar = \deg D - 1$. \\
(ii) The linear span of $\Dibar$ and $\Diibar$
is $\overline{\lcm(D_1,D_2)}$. \\
(iii) If $\deg (\lcm(D_1,D_2)) <n$ then $\Dibar \cap \Diibar =
\overline{\gcd(D_1,D_2)}$.
\end{Lemma}
\begin{Proof}
  (i) By Riemann-Roch we have $\dim \LL(H-D) = n - \deg D$. \\
  (ii) We have $\LL(H-D_1) \cap \LL(H-D_2) = \LL(H- \lcm(D_1,D_2))$. \\
  (iii) The inclusion ``$\supset$'' is clear. Equality follows
  by counting dimensions using~(i) and~(ii). 
\end{Proof}

With the above notation, Proposition~\ref{lem:tangent} becomes

\begin{Proposition}
\label{prop:tgt}
  Suppose $n-2r \ge 1$. Let $D = P_1 + \ldots + P_r$ be an effective
  divisor of degree $r$ with $P_1, \ldots,P_r \in C$ distinct. Then
  for any $P \in D^\circ$ we have $T_P \Sec^r C = \overline{2D}$.
\end{Proposition}

\begin{Proof}
  If $P \in \overline{D'}$ for $D'$ an effective divisor of degree at
  most $r$, then by Lemma~\ref{Dbarlemma}(iii) we have $D = D'$.  In
  particular $P \notin \Sec^{r-1}C$. It follows by
  Theorem~\ref{thm:hsec-eqns}(v) that $P$ is a smooth point on $\Sec^r
  C$.  The next lemma shows that $\overline{2D} \subset T_P \Sec^r C$,
  and equality follows by comparing dimensions, using
  Lemma~\ref{Dbarlemma}(i) and Theorem~\ref{thm:hsec-eqns}(i).
\end{Proof}

\begin{Lemma}
  Let $X$ be an affine variety and $P_1, \ldots, P_r \in X$.  Let $P =
  \sum \xi_i P_i$ where $\sum \xi_i =1$. If $\xi_i \not=0$ then
  $T_{P_i} X \subset T_P(\Sec^r X)$.
\end{Lemma}
\begin{Proof}
  There is a morphism $ X \times \ldots \times X \to \Sec^r X \,;\,
  (a_1, \ldots, a_r) \mapsto \sum \xi_i a_i$ with derivative $T_{P_1}X
  \times \ldots \times T_{P_r} X \to T_P(\Sec^r X) \,;\, (b_1, \ldots,
  b_r) \mapsto \sum \xi_i b_i$.
\end{Proof}

In fact Proposition~\ref{prop:tgt} is true without the hypothesis that
$P_1, \ldots, P_r$ are distinct. However, since we do not need this, 
we omit the details.

\subsection{Proof of Proposition~\ref{lem:codim3}}
\label{sec:codim3}

We must prove the following.

\begin{Proposition}
\label{codim3}
Suppose $n - 2 r = 2$ and write $\Sec^r C = \{F_1 = F_2 = 0\}$.  Then
the variety $X \subset \PP^{n-1}$ defined by
\[ \rank \begin{pmatrix} \vspace{0.5ex}
  \frac{\partial F_1}{\partial x_1} & \cdots & \frac{\partial F_1}{\partial x_n} \\
  \frac{\partial F_2}{\partial x_1} & \cdots & \frac{\partial
    F_2}{\partial x_n}
\end{pmatrix} \le 1 \]
has codimension $3$.
\end{Proposition}

If $n=4$ then $C = \{F_1 = F_2 = 0\} \subset \PP^3$ is the
intersection of two quadrics. There are $4$ singular quadrics in the
pencil spanned by $F_1$ and $F_2$, and each is singular at just one
point.  Then $X$ is the union of these $4$ singular points, and so has
codimension $3$.

We now generalise this argument. Let $H$ be the divisor of a
hyperplane section. We identify $\LL(H)$ with the space of linear
forms on $\PP^{n-1}$. Let $D_1$ and $D_2$ be divisors on $C$ of degree
$r+1$ with $D_1+D_2=H$. Let $\Phi(D_1,D_2)$ be the $(r+1) \times
(r+1)$ matrix of linear forms representing the multiplication map
\[ \LL(D_1) \times \LL(D_2) \to \LL(H). \] Since $\Phi(D_1,D_2)$ has
rank at most 1 on $C$, it has rank at most $r$ on $\Sec^r
C$. Therefore $\det \Phi(D_1,D_2)$ is a form of degree $r+1$ vanishing
on $\Sec^r C$. In particular it belongs to the pencil spanned by $F_1$
and $F_2$.

\begin{Lemma} 
\label{lem:every}
Every linear combination of $F_1$ and $F_2$ arises in this
way. Moreover there are exactly $4$ forms in the pencil arising as
$\det \Phi(D_1,D_2)$ with $D_1 \sim D_2$.
\end{Lemma}
\begin{Proof}
  We say that divisor pairs $(D_1,D_2)$ and $(D'_1,D'_2)$ are {\em
    equivalent} if $D_1 \sim D'_1$ or $D_1 \sim D'_2$. It is shown in
  \cite[Lemma~2.9]{pfpres} that if $(D_1,D_2)$ and $(D'_1,D'_2)$ are
  inequivalent then $\Sec^r C = \{ \det \Phi(D_1,D_2) = \det
  \Phi(D'_1,D'_2) = 0 \} \subset \PP^{n-1}$.  In particular these two
  forms are linearly independent.

  We claim that the map $(D_1,D_2) \mapsto \Phi(D_1,D_2)$ is a
  bijection between the equivalence classes of divisor pairs and the
  pencil of forms spanned by $F_1$ and $F_2$. To prove this let $C$ be
  the image of an elliptic curve $E$ embedded in $\PP^{n-1}$ by
  $|n.0_E|$. Then writing
  \[ \det \Phi(r.0_E + P,(r+2).0_E - P) = s(P) F_1 + t(P) F_2, \] for
  $P \in E$, we can see that $s/t$ is a rational function on $E$.  It
  therefore defines a morphism $(s:t) : E \to \PP^1$.  By the previous
  paragraph, this morphism is non-constant, and indeed has fibres of
  the form $\{P,-P\}$. It must therefore be surjective.  This proves
  the claim.

  For the final statement we note that $r.0_E + P \sim (r+2).0_E - P$
  if and only if $P \in E[2]$.
\end{Proof}

\begin{Lemma} 
  \label{lem:sing}
  Let $S$ be the singular locus of $V = \{ \det \Phi(D_1,D_2) = 0 \}
  \subset \PP^{n-1}$. Then $S$ contains $\Sec^{r-1} C$. Moreover
\begin{enumerate}
\item If $D_1 \not\sim D_2$ then $S  = \Sec^{r-1} C$. 
\item If $D_1 \sim D_2$ then $S$ has codimension $3$.
\end{enumerate}
\end{Lemma}

\begin{Proof}
  Since $C$ spans $\PP^{n-1}$ it is clear that for each $P \in
  \Sec^{r-1} C$ we have $T_P \Sec^r C = \PP^{n-1}$.
  Therefore $S$ contains $\Sec^{r-1} C$. \\
  (i) Let $P \in V \setminus \Sec^{r-1} C$ be any point.  According
  to~\cite[Theorem 1.3]{pfpres} the $r \times r$ minors of
  $\Phi(D_1,D_2)$ generate $I(\Sec^{r-1} C)$.  Therefore evaluating
  $\Phi(D_1,D_2)$ at $P$ gives a matrix of rank $r$.  Moving $P$ to
  $(1:0: \ldots:0)$ and picking suitable bases for $\LL(D_1)$ and
  $\LL(D_2)$ we have
  \[ \Phi(D_1,D_2) = x_1 \begin{pmatrix}
    0 & 0 & \cdots & 0 \\
    0 & 1 & \cdots & 0 \\
    \vdots & \vdots & \ddots & \vdots \\
    0 & 0 & \cdots & 1
  \end{pmatrix} + \Phi' \] where $\Phi'$ is an $(r+1) \times (r+1)$
  matrix of linear forms in $x_2, \ldots, x_n$.  Now the top left
  entry of $\Phi(D_1,D_2)$ is an equation for $T_P V$.  Since the
  product of non-zero rational functions on $C$ is again non-zero, the
  entries of $\Phi(D_1,D_2)$ are non-zero. Therefore
  $P \in V$ is a smooth point. \\
  (ii) Picking suitable bases for $\LL(D_1)$ and $\LL(D_2)$ we may
  suppose that $\Phi(D_1,D_2)$ is symmetric. Since $\{ \rank
  \Phi(D_1,D_2) \le r-1 \} \subset S$, and the quadratic forms of rank
  at most~$m-2$ have codimension~$3$ in the space of all quadratic
  forms in $m$ variables, it follows that $S$ has codimension at
  most~$3$.  Suppose for a contradiction that $S$ has codimension at
  most~$2$. Then its intersection with $\Sec^r C = \{F_1 = F_2 = 0 \}$
  has codimension at most~$3$. But this intersection is contained in
  the singular locus of $\Sec^r C$, which by
  Theorem~\ref{thm:hsec-eqns} has codimension~$4$.  This is the
  required contradiction.
\end{Proof}

To complete the proof of Proposition~\ref{codim3}, we note that $X$
is the union of the singular loci of the hypersurfaces defined 
by linear combinations of $F_1$ and $F_2$. It follows by 
Lemmas~\ref{lem:every} and~\ref{lem:sing} that $X$ has codimension $3$.


\begin{thebibliography}{MM}

\frenchspacing
\renewcommand{\baselinestretch}{1}

\bibitem[AKM$^3$P]{Mc+}
S.Y. An, S.Y. Kim, D.C. Marshall, S.H. Marshall, W.G. McCallum and 
A.R. Perlis, 
Jacobians of genus one curves,
{\em J. Number Theory} 90 (2001), no. 2, 304--315.

\bibitem[A]{Aronhold}
S. Aronhold, Theorie der homogenen Funktionen dritten Grades 
von drei Ver\"anderlichen, 
{\em J. reine angew. Math.} 55 (1858), 97--191.

\bibitem[ARVT]{ARVT}
M. Artin, F. Rodriguez-Villegas and J. Tate, 
On the Jacobians of plane cubics,
{\em Adv. Math.} 198 (2005), no. 1, 366--382.

\bibitem[B]{quintic}
M. Bhargava, 
Higher composition laws, IV. The parametrization of quintic rings,
{\em Ann. of Math.} (2) 167 (2008), no. 1, 53--94.

\bibitem[vBH]{vBH}
H.-Chr. Graf v. Bothmer and K. Hulek, 
Geometric syzygies of elliptic normal curves and their secant varieties,
{\em Manuscripta Math.} {{113}} (2004), no. 1, 35--68.

\bibitem[BH]{CMrings}
W. Bruns and J. Herzog, 
{\em Cohen-Macaulay rings},
Cambridge Studies in Advanced Mathematics, 39, 
Cambridge University Press, Cambridge, 1993. 

\bibitem[BE1]{BE1}
D.A. Buchsbaum and D. Eisenbud, 
Gorenstein ideals of height $3$,
{\em Seminar D. Eisenbud/B. Singh/W. Vogel}, Vol. 2, 
pp. 30--48, Teubner-Texte zur Math., 48, Teubner, Leipzig, 1982. 

\bibitem[BE2]{BE2}
D.A. Buchsbaum and D. Eisenbud,
Algebra structures for finite free resolutions, and some 
structure theorems for ideals of codimension 3,
{\em Amer. J. Math.} 
{{99}} (1977) 447-485.

\bibitem[E]{Eisenbud}
D. Eisenbud, 
{\em Commutative algebra with a view toward algebraic geometry},
Graduate Texts in Mathematics, 150, Springer-Verlag, New York, 1995. 

\bibitem[F1]{g1inv}
T.A. Fisher,
The invariants of a genus one curve, 
{\em Proc. Lond. Math. Soc.} (3) {97} (2008) 753--782. 

\bibitem[F2]{pfpres}
T.A. Fisher, 
Pfaffian presentations of elliptic normal curves,
{\em Trans. Amer. Math. Soc.} 362 (2010), no. 5, 2525--2540.

\bibitem[F3]{invenqI}
T.A. Fisher,
Invariant theory for the elliptic normal quintic, I. Twists of $X(5)$,
{\em Math. Ann.} 356 (2013), no. 2, 589--616.

\bibitem[F4]{5desc}
T.A. Fisher,
Explicit 5-descent on elliptic curves, in
{\em ANTS X -- Proceedings of the Tenth Algorithmic Number 
Theory Symposium}, E.W. Howe and K.S. Kedlaya (eds), 
Open Book Ser., 1, Math. Sci. Publ., Berkeley, CA, 2013, 395--411.

\bibitem[FS]{sqrfree}
T.A. Fisher and M. Sadek,
{\em On genus one curves of degree 5 with square-free discriminant}, 
preprint, 2014.
 
\bibitem[G]{BGross}
B.H. Gross, 
On Bhargava's representation and Vinberg's invariant theory,
in {\em Frontiers of mathematical sciences}, 
B. Gu and S.-T. Yau (eds), Int. Press, Somerville, MA, 2011, 317--321.

\bibitem[GP]{GP}
M. Gross and S. Popescu, 
Equations of $(1,d)$-polarized abelian surfaces, 
{\em Math. Ann.} {{310}} (1998), no. 2, 333--377.

\bibitem[H]{Hulek} 
K. Hulek,
{\em Projective geometry of elliptic curves},
Soc. Math. de France, Ast\'erisque
{{137}} (1986).

\bibitem[L]{Lange}
H. Lange, 
Higher secant varieties of curves and the theorem of Nagata on ruled surfaces,
{\em Manuscripta Math.} {{47}} (1984), no. 1-3, 263--269. 

\bibitem[P]{Peeva}
I. Peeva, 
{\em Graded syzygies},
Algebra and Applications, 14, Springer-Verlag, London, 2011. 

\bibitem[S]{Sil}
J.H. Silverman, 
{\em The arithmetic of elliptic curves}, Graduate Texts 
in Mathematics, 106, Springer, Dordrecht, 2009. 

\bibitem[W1]{WHermite}
A. Weil, 
Remarques sur un m\'emoire d'Hermite,
{\em Arch. Math.} (Basel) 5, (1954), 197--202.

\bibitem[W2]{WEuler}
A. Weil, 
Euler and the Jacobians of elliptic curves, {\em Arithmetic and geometry}, 
Vol. I, 353--359, Progr. Math., 35, Birkh\"auser, Boston, 1983.

\end{thebibliography}
\end{document}